\documentclass[11pt]{article}
\usepackage{amssymb}
\usepackage{amsfonts}
\usepackage{amsmath}
\usepackage{bm}
\usepackage{latexsym}
\usepackage{epsfig}

\newtheorem{theorem}{Theorem}[section]
\newtheorem{lemma}[theorem]{Lemma}
\newtheorem{proposition}[theorem]{Proposition}
\newtheorem{corollary}[theorem]{Corollary}

\newtheorem{definition}[theorem]{Definition}
\newtheorem{remark}[theorem]{Remark}
\newtheorem{conjecture}[theorem]{Conjecture}
\newenvironment{proof}{\noindent \textbf{Proof:}}{$\Box$}

\newcommand{\ignore}[1]{}

\newcommand{\enote}[1]{}
\newcommand{\knote}[1]{}
\newcommand{\rnote}[1]{}

\renewcommand{\Pr}{{\bf P}}
\renewcommand{\P}{{\bf P}}
\newcommand{\Px}{\mathop{\bf P\/}}
\newcommand{\E}{{\bf E}}
\newcommand{\Ex}{\mathop{\bf E\/}}
\newcommand{\Var}{{\bf Var}}
\newcommand{\Varx}{\mathop{\bf Var\/}}

\newcommand{\bits}{\{-1,1\}}

\newcommand{\Inf}{\mathrm{Inf}}

\newcommand{\eps}{\epsilon}
\newcommand{\lam}{\lambda}

\newcommand{\trunc}{\zeta}

\newcommand{\N}{\mathbb N}
\newcommand{\R}{\mathbb R}
\newcommand{\CalE}{{\mathcal{E}}}
\newcommand{\CalC}{{\mathcal{C}}}

\newcommand{\CalX}{{\boldsymbol{\mathcal{X}}}}
\newcommand{\CalG}{{\boldsymbol{\mathcal{G}}}}
\newcommand{\CalY}{{\boldsymbol{\mathcal{Y}}}}
\newcommand{\CalZ}{{\boldsymbol{\mathcal{Z}}}}

\newcommand{\boldG}{{\boldsymbol G}}
\newcommand{\boldQ}{{\boldsymbol Q}}
\newcommand{\boldR}{{\boldsymbol R}}
\newcommand{\boldS}{{\boldsymbol S}}
\newcommand{\boldX}{{\boldsymbol X}}
\newcommand{\boldB}{{\boldsymbol B}}
\newcommand{\boldY}{{\boldsymbol Y}}
\newcommand{\boldZ}{{\boldsymbol Z}}
\newcommand{\boldV}{{\boldsymbol V}}
\newcommand{\boldsigma}{{\boldsymbol \sigma}}
\newcommand{\boldupsilon}{{\boldsymbol \upsilon}}
\newcommand{\hone}{{\boldsymbol{H1}}}
\newcommand{\htwo}{\boldsymbol{H2}}
\newcommand{\hthree}{\boldsymbol{H3}}
\newcommand{\hfour}{\boldsymbol{H4}}

\newcommand{\Maj}{\mathrm{Maj}}
\newcommand{\Thr}{\mathrm{Thr}}
\newcommand{\littlesum}{{\textstyle \sum}}

\newcommand{\half}{{\textstyle \frac12}}

\newcommand{\Stab}{\mathbb{S}}
\newcommand{\StabThr}[2]{\Gamma_{#1}(#2)}
\newcommand{\TestFcn}{\Psi}

\renewcommand{\phi}{\varphi}

\begin{document}
\title{Noise stability of functions with low influences: \linebreak invariance
  and optimality }

\author{Elchanan Mossel\footnote{Supported by a Miller fellowship in
Statistics and Computer Science, U.C. Berkeley and by a Sloan
fellowship in Mathematics}
\\U.C.\ Berkeley\\mossel@stat.berkeley.edu \and Ryan
O'Donnell\\Microsoft Research\\odonnell@microsoft.com \and
Krzysztof Oleszkiewicz\footnote{Partially supported by  Polish KBN
Grant 2 P03A 027 22.}\\Warsaw University\\koles@mimuw.edu.pl}
\date{\today}
\maketitle

\begin{abstract}
In this paper we study functions with low influences on
product probability spaces.  These are functions $f : \Omega_1
\times \cdots \times \Omega_n \to \R$ that have
$\E[\Var_{\Omega_i}[f]]$ small compared to $\Var[f]$ for each
$i$. The analysis of boolean functions $f: \{-1,1\}^n \to
\{-1,1\}$ with low influences has become a central problem in
discrete Fourier analysis. It is motivated by fundamental
questions arising from the construction of probabilistically
checkable proofs in theoretical computer science and from problems
in the theory of social choice in economics.

We prove an invariance principle for multilinear polynomials with
low influences and bounded degree; it shows that under mild conditions 
the distribution of such polynomials is essentially invariant 
for all product spaces.  Ours is one of the very few known non-linear
invariance principles. It has the advantage that its proof is
simple and that the error bounds are explicit. We also show that
the assumption of bounded degree can be eliminated if the
polynomials are slightly ``smoothed''; this extension is essential
for our applications to ``noise stability''-type problems.

In particular, as applications of the invariance principle we
prove two conjectures: the ``Majority Is
Stablest'' conjecture~\cite{KKMO:04} from theoretical computer
science, which was the original motivation for this work, and the
``It Ain't Over Till It's Over'' conjecture~\cite{Kalai:01} from
social choice theory.
\end{abstract}

\section{Introduction}

\subsection{Harmonic analysis of boolean functions}
\label{sec:intro}  The motivation for this paper is the study of
\emph{boolean functions} $f : \bits^n \to \bits$, where $\bits^n$
is equipped with the uniform probability measure. This topic is of
significant interest in theoretical computer science; it also
arises in other diverse areas of mathematics including
combinatorics (e.g., sizes of set systems, additive
combinatorics), economics (e.g., social choice), metric spaces
(e.g., non-embeddability of metrics), geometry in Gaussian space
(e.g., isoperimetric inequalities), and
statistical physics (e.g., percolation, spin glasses).\\

Beginning with Kahn, Kalai, and Linial's landmark paper ``The
Influence Of Variables On Boolean Functions''~\cite{KaKaLi:88}
there has been much success in analyzing questions about boolean
functions using methods of harmonic analysis. Recall that KKL
essentially shows the following (see
also~\cite{Talagrand:94,FriedgutKalai:96}):

\paragraph{KKL Theorem:} If $f : \bits^n
\to \bits$ satisfies $\E[f] = 0$ and $\Inf_i(f) \leq \tau$ for all
$i$, then $\sum_{i=1}^n \Inf_i(f) \geq \Omega(\log(1/\tau))$.\\

\noindent We have used here the notation $\Inf_i(f)$ for the
\emph{influence of the $i$th coordinate on $f$},
\begin{equation} \label{eqn:influence}
\Inf_i(f) = \Ex_x[\Varx_{x_i}[f(x)]] = \sum_{S \ni i}
\hat{f}(S)^2.
\end{equation}

Although an intuitive understanding of the analytic properties of
boolean functions is emerging, results in this area have used
increasingly elaborate methods, combining random restriction
arguments, applications of the Bonami-Beckner inequality, and
classical tools from probability theory. See for
example~\cite{Talagrand:94,Talagrand:96,FriedgutKalai:96,Friedgut:99,Bourgain:99,BeKaSc:99,Bourgain:02,KindlerSafra:u,DiFrKiOD:u}.\\

As in the KKL paper, some of the more refined problems studied in
recent years have involved restricting attention to functions with
low influences~\cite{BeKaSc:99,Bourgain:99,DiFrKiOD:u} (or,
relatedly, ``non-juntas''). There are several reasons for this.
The first is that large-influence functions such as ``dictators''
--- i.e., functions $f(x_1, \dots, x_n) = \pm x_i$ --- frequently trivially
maximize or minimize quantities studied in boolean analysis.
However this tends to obscure the truth about extremal behaviors
among functions that are ``genuinely'' functions of $n$ bits.
Another reason for analyzing only low-influence functions is that
this subclass is often precisely what is interesting or necessary
for applications.  In particular, the analysis of low-influence
boolean functions is crucial for proving hardness of approximation
results in theoretical computer science and is also very natural
for the study of social choice.  Let us describe these two
settings briefly.\\

In the economic theory of social choice, boolean functions $f :
\bits^n \to \bits$ often represent voting schemes, mapping $n$
votes between two candidates into a winner.  In this case, it is
very natural to exclude voting schemes that give any voter an
undue amount of influence; see e.g.~\cite{Kalai:04}. In the study
of hardness of approximation and probabilistically checkable
proofs (PCPs), the sharpest results often involve the following
paradigm: One considers a problem that requires labeling the
vertices of a graph using the label set $[n]$; then one relaxes
this to the problem of labeling the vertices by functions $f :
\bits^n \to \bits$. In the relaxation one thinks of $f$ as
``weakly labeling'' a vertex by the \emph{set} of coordinates that
have large influence on $f$.  It then becomes important to
understand the combinatorial properties of functions that weakly
label with the empty set.  There are by now quite a few results in
hardness of approximation that use results on low-influence
functions or require conjectured such results; e.g.,
\cite{DinurSafra:02,Khot:02,DGKO:03,KhotRegev:03,KKMO:04}.\\

In this paper we give a new framework for studying functions on
product probability spaces with low influences.  Our main tool is
a simple invariance principle for low-influence polynomials; this
theorem lets us take an optimization problem for functions on one
product space and pass freely to other product spaces, such as
Gaussian space.  In these other settings the problem sometimes
becomes simpler to solve.  It is interesting to note that while in
the theory of hypercontractivity and isoperimetry it is common to
prove results in the Gaussian setting by proving them first in the
$\bits^n$ setting (see, e.g., \cite{Bakry:94}),
here the invariance principle is actually used to go the other way around.\\

As applications of our invariance principle we prove two
previously unconnected conjectures from boolean harmonic analysis;
the first was motivated by hardness of approximation in computer
science, the second by the theory of social choice from economics:

\begin{conjecture}[``Majority Is Stablest''
  conjecture~\cite{KKMO:04}] \label{conj:MIST}
Let $0 \leq \rho \leq 1$ and $\eps > 0$ be given.  Then there
exists $\tau > 0$ such that if $f : \bits^n \to [-1,1]$ satisfies
$\E[f] = 0$ and $\Inf_i(f) \leq \tau$ for all $i$, then
\[
\Stab_\rho(f) \leq {\textstyle \frac{2}{\pi}} \arcsin \rho + \eps.
\]
\end{conjecture}
Here we have used the notation $\Stab_\rho(f)$ for $\sum_S
\rho^{|S|} \hat{f}(S)^2$, the \emph{noise stability} of $f$.  This
equals $\E[f(x)f(y)]$ when $(x,y) \in \bits^n \times \bits^n$
is chosen so that $(x_i,y_i) \in \bits^2$ are independent random
variables with $\E[x_i] = \E[y_i] = 0$ and $\E[x_i y_i] = \rho$.

\begin{conjecture}[``It Ain't Over Till It's Over''
conjecture~\cite{Kalai:01}] \label{conj:ain't} Let $0 \leq \rho <
1$ and $\eps > 0$ be given.  Then there exists $\delta > 0$ and
$\tau > 0$ such that if $f : \bits^n \to \bits$ satisfies $\E[f] =
0$ and $\Inf_i(f) \leq \tau$ for all $i$, then $f$ has the
following property:  If $V$ is a random subset of $[n]$ in which
each $i$ is included independently with probability $\rho$, and if
the bits $(x_i)_{i \in V}$ are chosen uniformly at random, then
\[
\Px_{V,\;(x_i)_{i \in V}} \Bigl[ \bigl|\E[f \mid (x_i)_{i \in
V}]\bigr| > 1 - \delta\Bigr] \leq \eps.
\]
\end{conjecture}
(In words, the conjecture states that even if a random $\rho$
fraction of voters' votes are revealed, with high probability the
election is still slightly undecided, provided $f$ has low
influences.)\\

The truth of these results gives illustration to a recurring theme
in the harmonic analysis of boolean functions: the extremal role
played the Majority function. It seems this theme becomes
especially prominent when low-influence functions are studied. To
explain the connection of Majority to our applications: In the
former case the quantity $\frac{2}{\pi} \arcsin \rho$ is precisely
$\lim_{n \to \infty} \Stab_\rho(\Maj_n)$; this explains the name
of the Majority Is Stablest conjecture. In the latter case, we
show that $\delta$ can be taken to be on the order of
$\epsilon^{\rho/(1-\rho)}$ (up to $o(1)$ in the exponent), which
is the same asymptotics one gets if $f$ is Majority on a large
number of inputs.

\subsection{Outline of the paper}
We begin in Section~\ref{sec:statement} with an overview of the
invariance principle, the two applications, and some of their consequences.  
We prove the invariance principle in
Section~\ref{sec:invariance}. Our proofs of the two conjectures
are in Section~\ref{sec:conj}. Finally, we show in
Section~\ref{sec:counterexample} that a conjecture closely related
to Majority Is Stablest is false.  Some minor proofs from
throughout the paper appear in appendices.

\subsection{Related work}
Our multilinear invariance principle has some antecedents.  For
degree 1 polynomials it reduces to a version of the
Berry-Esseen Central Limit Theorems.  Indeed, our proof
follows the same outlines as Lindeberg's proof of the
CLT~\cite{Lindeberg:22} (see also~\cite{Feller:71}).\\

Since presenting our proof of the invariance principle, we have
been informed by Oded Regev that related results were proved in
the past by V.~I.~Rotar$'$~\cite{Rotar:79}.  As well, a
contemporary manuscript of Sourav Chatterjee~\cite{Chatterjee:u}
with an invariance principle of similar flavor has come to our
attention. What is common to our work and
to~\cite{Rotar:79,Chatterjee:u} is a generalization of Lindeberg's
argument to the non-linear case. The result of Rotar$'$ is an
invariance principle similar to ours where the condition on the
influences generalizes Lindeberg's condition.  The setup is not
quite the same, however, and the proof in~\cite{Rotar:79} is of a
rather qualitative nature. It seems that even after appropriate
modification the bounds it gives would be weaker and less useful
for our type of applications.  (This is quite understandable; in a
similar way Lindeberg's CLT can be less precise than the
Berry-Esseen inequality even though --- indeed, because --- it
works under weaker assumptions.)  The paper~\cite{Chatterjee:u} is
by contrast very clear and explicit. However it does not seem to
be appropriate for many applications since it requires low
``worst-case'' influences, instead of the ``average-case''
influences used by this work and~\cite{Rotar:79}.\\

Finally, we would like to mention that some chaos-decomposition
limit theorems have been proved before in various settings.  Among
these are limit theorems for U and V statistics and limit theorems
for random graphs; see, e.g.~\cite{Janson:97}.

\subsection{Acknowledgments}
We are grateful to Keith Ball for suggesting a collaboration among
the authors.  We would also like to thank Oded Regev for referring
us to~\cite{Rotar:79} and Olivier Gu\'edon for referring us
to~\cite{CarberyWright:01}.

\section{Our results} \label{sec:statement}

\subsection{The invariance principle} \label{sec:intro-invariance}
In this subsection we present a simplified version of our invariance principle.\\

Suppose $\boldX$ is a random variable with $\E[\boldX] = 0$ and
$\E[\boldX^2] = 1$ and $\boldX_1, \dots, \boldX_n$ are independent
copies of $\boldX$.  Let $Q(x_1, \dots, x_n) = \sum_{i=1}^n c_i
x_i$ be a linear form and assume $\sum c_i^2 = 1$.  The
Berry-Esseen CLT states that under mild conditions on the
distribution of $\boldX$, say $\E[|\boldX|^3] \leq A < \infty$, it
holds that
\[
\sup_{t} \bigl| \P[Q(\boldX_1, \dots, \boldX_n) \leq t] -
          \P[\boldG \leq t]\bigr| \leq O\bigl(A \cdot \littlesum_{i=1}^n |c_i|^3\bigr),
\]
where $\boldG$ denotes a standard normal random variable.  Note
that a simple corollary of the above is
\begin{equation}  \label{eq:berry}
\sup_{t} \bigl|\P[Q(\boldX_1, \dots, \boldX_n) \leq t] -
          \P[Q(\boldG_1, \dots, \boldG_n) \leq t] \bigr| \leq
          O\bigl(A \cdot \max_i |c_i|\bigr).
\end{equation}
Here the $\boldG_i$'s denote independent standard normals.  We
have upper-bounded the sum of $|c_i|^3$ here by a maximum, for
simplicity; more importantly though, we have suggestively replaced
$\boldG$ by $\sum_i c_i \boldG_i$, which of course has the same
distribution.\\

We would like to generalize~(\ref{eq:berry}) to \emph{multilinear
polynomials} in the $\boldX_i$'s; i.e., functions of the form
\begin{equation} \label{eqn:Q1}
Q(\boldX_1, \dots, \boldX_n) = \sum_{S \subseteq [n]} c_S \prod_{i
\in S} \boldX_i,
\end{equation}
where the real constants $c_S$ satisfy $\sum c_S^2 = 1$.  Let $d =
\max_{c_S \neq 0} |S|$ denote the degree of $Q$.  Unlike in the
$d=1$ case of the CLT, there is no single random variable $\boldG$
which always provides a limiting distribution. However one can
still hope to prove, in light of~(\ref{eq:berry}), that the
distribution of the polynomial applied to the variables $\boldX_i$
is close to the distribution of the polynomial applied to
independent Gaussian random variables. This is indeed what our
invariance principle shows.\\

It turns out that the appropriate generalization of the
Berry-Esseen theorem~(\ref{eq:berry}) is to control the error
by a function of $d$ and of $\max_i \sum_{S \ni i} c_S^2$ ---
i.e., the maximum of the \emph{influences} of $Q$ (as
in~(\ref{eqn:influence})).  Naturally, we also
need some conditions in addition to second moments.  In our
formulation we impose the condition that the variable $\boldX$ is
\emph{hypercontractive}; i.e., there is some $\eta > 0$ such that
for all $a \in \R$,
\[
\|a + \eta \boldX\|_3 \leq \|a + \boldX\|_2.
\]
This condition is satisfied whenever $\E[\boldX]=0$ and
$\E[|\boldX|^{3}]<\infty;$ in particular, it holds for any
mean-zero random variable $\boldX$ taking on only finitely many
values. Using hypercontractivity, we get a simply proved
invariance principle with explicit error bounds.  The following
theorem (a simplification of Theorem~\ref{thm:supertheorem},
bound~(\ref{eq:lim_dist})) is an example of what we prove:
\begin{theorem} \label{thm:simple}
Let $\boldX_1, \dots, \boldX_n$ be independent random variables
satisfying $\E[\boldX_i] = 0$, $\E[\boldX_i^2] = 1$, and
$\E[|\boldX_i|^{3}] \leq \beta.$ Let $Q$ be a degree $d$
multilinear polynomial as in~(\ref{eqn:Q1}) with
\[
\sum_{|S| > 0} c_S^2 = 1, \qquad \qquad  \sum_{S \ni i} c_S^2 \leq
\tau \quad \text{for all $i$}.
\]
Then
\[
\sup_{t} \bigl|\P[Q(\boldX_1, \dots, \boldX_n) \leq t] -
          \P[Q(\boldG_1, \dots, \boldG_n) \leq t] \bigr|
  \leq O(d\beta^{1/3}\tau^{1/8d}),
\]
where $\boldG_1, \dots, \boldG_n$ are independent standard
Gaussians.\\

If, instead of assuming $\E[|\boldX_i|^{3}] \leq \beta$, we assume
that each $\boldX_i$ takes only on finitely many values, and that
for all $i$ and all $x \in \R$ either $\Pr[\boldX_i = x] = 0$ or
$\Pr[\boldX_i = x] \geq \alpha$, then
\[
\sup_{t} \bigl|\P[Q(\boldX_1, \dots, \boldX_n) \leq t] -
          \P[Q(\boldG_1, \dots, \boldG_n) \leq t] \bigr|
  \leq O(d\,\alpha^{-1/6}\, \tau^{1/8d}).
\]
\end{theorem}
Note that if $d$, $\beta$, and $\alpha$ are fixed then the above
bound tends to $0$ with $\tau$.  We call this theorem an
``invariance principle'' because it shows that $Q(\boldX_1, \dots,
\boldX_n)$ has essentially the same distribution no matter what
the $\boldX_i$'s are. Usually we will not push for the optimal
constants; instead we will try to keep our approach as simple as
possible while still giving explicit bounds useful for our
applications.\\

An unavoidable deficiency of this sort of invariance principle is
the dependence on $d$ in the error bound.  In applications such as
Majority Is Stablest and It Ain't Over Till It's Over, the
functions $f$ may well have arbitrarily large degree.  To overcome
this, we introduce a supplement to the invariance principle:  We
show that if the polynomial $Q$ is ``smoothed'' slightly then the
dependence on $d$ in the error bound can be eliminated and
replaced with a dependence on the smoothness.  For ``noise
stability''-type problems such as ours, this smoothing is
essentially harmless.\\

In fact, the techniques we use are strong enough to obtain Berry-Esseen 
estimates under Lyapunov-type assumptions. In particular, we believe that the 
following theorem is new even in the case of sums of independent 
random variables. 

\begin{theorem} \label{thm:Lyap}
Let $q \in (2,3].$
Let $\boldX_1, \dots, \boldX_n$ be independent random variables
satisfying $\E[\boldX_i] = 0$, $\E[\boldX_i^2] = 1$, and
$\E[|\boldX_i|^{q}] \leq \beta.$ Let $Q$ be a degree $d$
multilinear polynomial as in~(\ref{eqn:Q1}) with
\[
\sum_{|S| > 0} c_S^2 = 1, \qquad \qquad  \sum_{S \ni i} c_S^2 \leq
\tau \quad \text{for all $i$}.
\]
Then
\[
\sup_{t} \bigl|\P[Q(\boldX_1, \dots, \boldX_n) \leq t] -
          \P[Q(\boldG_1, \dots, \boldG_n) \leq t] \bigr|
  \leq 
  \]
  \[
  O(d\beta^{\frac{d}{qd+1}}
  (\sum_{i}(\sum_{S \ni i} c_{S}^{2})^{q/2})^{\frac{1}{qd+1}})
  \leq
   O(d\beta^{\frac{d}{qd+1}}
  \tau^{\frac{q-2}{2qd+2}}),
\]
where $\boldG_1, \dots, \boldG_n$ are independent standard
Gaussians.\\
\end{theorem}

\subsection{Influences and noise stability in product
spaces} \label{sec:general}

Our proofs of the Majority Is Stablest and It Ain't Over Till It's
Over conjectures hold not just for functions on the
uniform-distribution discrete cube, but for functions on arbitrary
finite product probability spaces.  Harmonic analysis results on
influences have often considered the biased product distribution
on the discrete cube (see, e.g.,
\cite{Talagrand:94,FriedgutKalai:96,Friedgut:99,Bourgain:99});
and, some recent works involving influences and noise stability
have considered functions on product sets $[q]^n$ endowed with the
uniform distribution (e.g., \cite{AlDiFrSu:04,KKMO:04}).  But since
there doesn't
appear to be a unified treatment for the general case in the literature, we give the necessary definitions here.\\

Let $(\Omega_1, \mu_1), \dots, (\Omega_n, \mu_n)$ be
probability spaces and let $(\Omega, \mu)$ denote the
product probability space.  Let
\[
f : \Omega_1 \times \cdots \times \Omega_n \to \R
\]
be any real-valued function on $\Omega$.

\begin{definition} \label{def:influence_general}
The \emph{influence of the $i$th coordinate on $f$} is
\[
\Inf_i(f) = \Ex_{\mu} [ \Varx_{\mu_i} [f]].
\]
\end{definition}
Note that for boolean functions $f : \bits^n \to \bits$ this
agrees with the classical notion of influences introduced to
computer science by Ben-Or and Linial~\cite{BenorLinial:90}.  When
the domain $\bits^n$ has a $p$-biased distribution, our notion
differs from that of, say,~\cite{Friedgut:98} by a multiplicative
factor of $4p(1-p)$.  We believe the above definition is more
natural, and in any case it is easy to pass between the two.\\

To define noise stability, we first define the $T_\rho$ operator
on the space of functions $f$:
\begin{definition} \label{def:T_general}
For any $0 \leq \rho \leq 1$, the operator $T_\rho$ is defined by
\begin{equation} \label{eq:T_general}
(T_\rho f)(\omega_1, \dots, \omega_n) = \E[f(\omega_1', \dots,
\omega_n')],
\end{equation}
where each $\omega_i'$ is an independent random variable defined
to equal $\omega_i$ with probability $\rho$ and to be randomly
drawn from $\mu_i$ with probability $1-\rho$.
\end{definition}

\medskip
We remark that this definition agrees with that of the
``Bonami-Beckner operator'' introduced in the context of boolean
functions by KKL~\cite{KaKaLi:88} and also with its generalization
to $[q]^n$ from~\cite{KKMO:04}.  For more on this operator, see
Wolff~\cite{Wolff:u}.  With this definition in place, we can
define noise stability:
\begin{definition} \label{def:Stab_general}
The \emph{noise stability of $f$ at $\rho \in [0,1]$} is
\[
\Stab_\rho(f) = \Ex_\mu[f \cdot T_\rho f].
\]
\end{definition}

\bigskip

For the It Ain't Over Till It's Over problem, we introduce a new
operator $V_\rho$:
\begin{definition} \label{def:V}
For any $\rho \in [0,1]$, the operator $V_\rho$ is defined as
follows. The operator takes a function 
$f : \Omega_1 \times \cdots \times \Omega_n \to \R$ to a function 
$g : \Omega_1 \times \cdots \times \Omega_n \times \{0,1\}^n \to \R$, where
$\{0,1\}^n$ is equipped with the $(1-\rho,\rho)^{\otimes n}$ measure. 
It is defined as follows: 
\[
(V_\rho f)(\omega_1, \dots, \omega_n,x_1,\ldots,x_n) = 
\Ex_{\omega'} \left[
f \left( x_1 \omega_1 + (1-x_1) \omega'_1, \ldots, 
x_n \omega_n + (1-x_n) \omega'_n \right) \right].
\]
\end{definition}

\medskip
Finally, we would like to note that our definitions are valid for
functions $f$ into the reals, although our motivation is usually
$\bits$-valued functions.  Our proofs of the Majority Is Stablest
and It Ain't Over Till It's Over conjectures will hold in the
setting of functions $f : \Omega_1 \times \cdots \times \Omega_n
\to [-1,1]$ (note that Conjecture~\ref{conj:MIST} \emph{requires}
this generalized range).  For notational simplicity, though, we
will give our proofs for functions into $[0,1]$; the reader can
easily convert such results to the $[-1,1]$ case by the linear
transformation $f \mapsto 2f-1$, which interacts in a simple way
with the definitions of $\Inf_i$, $\Stab_\rho$ and $V_\rho$.

\subsection{Majority Is Stablest} \label{sec:misc}
\subsubsection{About the problem} \label{sec:misc-discuss}
The Majority Is Stablest conjecture, Conjecture~\ref{conj:MIST},
was first formally stated in~\cite{KKMO:04}. However the notion of
Hamming balls having the highest noise stability in various senses
has been widely spread among the community studying discrete
Fourier analysis.  Indeed, already in KKL's 1998
paper~\cite{KaKaLi:88} there is the suggestion that Hamming balls
and subcubes should maximize a certain noise stability-like
quantity. In~\cite{BeKaSc:99}, it was shown that every
`asymptotically noise stable'' function is correlated with a
weighted majority function; also, in~\cite{MORSS:04} it was shown
that the majority function asymptotically maximizes a high-norm
analog of $\Stab_{\rho}$.\\

More concretely, strong motivation for getting
sharp bounds on the noise stability of low-influence functions
came from two 2002 papers, one by Kalai~\cite{Kalai:02} on social
choice and one by Khot~\cite{Khot:02} on PCPs and hardness of
approximation.  We briefly discuss these two papers below.\\

\paragraph{Kalai '02 --- Arrow's Impossibility Theorem:} Suppose $n$ voters
rank three candidates, $A$, $B$, and $C$, and a \emph{social
choice} function $f : \bits^n \to \bits$ is used to aggregate the
rankings, as follows: $f$ is applied to the $n$ $A$-vs.-$B$
preferences to determine whether $A$ or $B$ is globally preferred;
then the same happens for $A$-vs.-$C$ and $B$-vs.-$C$.  The
outcome is termed ``non-rational'' if the global ranking has $A$
preferable to $B$ preferable to $C$ preferable to $A$ (or if the
other cyclic possibility occurs).  Arrow's Impossibility Theorem
from the theory of social choice states that under some mild
restrictions on $f$ (such as $f$ being odd; i.e., $f(-x) =
-f(x)$), the only functions which never admit non-rational
outcomes given rational voters are the dictator functions $f(x) =
\pm x_i$.

Kalai~\cite{Kalai:02} studied the \emph{probability} of a rational
outcome given that the $n$ voters vote independently and at random
from the 6 possible rational rankings.  He showed that the
probability of a rational outcome in this case is precisely $3/4 +
(3/4) \Stab_{1/3}(f)$.  Thus it is natural to ask which function
$f$ with small influences is most likely to produce a rational
outcome.  Instead of considering small influences, Kalai
considered the essentially stronger assumption that $f$ is
``transitive-symmetric''; i.e., that for all $1 \leq i < j \leq n$
there exists a permutation $\sigma$ on $[n]$ with $\sigma(i) = j$
such that $f(x_1, \dots, x_n) = f(x_{\sigma(1)}, \dots,
x_{\sigma(n)})$ for all $(x_1, \dots, x_n)$.  Kalai conjectured
that Majority was the transitive-symmetric function that maximized
$3/4 + (3/4) \Stab_{1/3}(f)$ (in fact, he made a stronger
conjecture, but this conjecture is false; see
Section~\ref{sec:counterexample}). He further observed that this
would imply that in any transitive-symmetric scheme the
probability of a rational outcome is at most $3/4 +
(3/2\pi)\arcsin(1/3) + o_n(1) \approx .9123$; however, Kalai could
only prove the weaker bound $.9192$.

\paragraph{Khot '02 --- Unique Games and hardness of approximating
2-CSPs:}  In computer science, many combinatorial optimization
problems are NP-hard, meaning it is unlikely there are efficient
algorithms that always find the optimal solution.  Hence there has
been extensive interest in understanding the complexity of
\emph{approximating} the optimal solution.  Consider for example
``$k$-variable constraint satisfaction problems'' ($k$-CSPs) in
which the input is a set of variables over a finite domain,
along with some constraints on $k$-sets of the variables,
restricting what sets of values they can simultaneously take.  We
say a problem has ``$(c,s)$-hardness'' if it is NP-hard, given a
$k$-CSP instance in which the optimal assignment satisfies a
$c$-fraction of the constrains, for an algorithm to find an
assignment that satisfies an $s$-fraction of the constraints.  In
this case we also say that the problem is ``$s/c$-hard to
approximate''.

The PCP and Parallel Repetition theorems have led to many
impressive results showing that it is NP-hard even to give
$\alpha$-approximations for various problems, especially $k$-CSPs
for $k \geq 3$. For example, letting MAX-$k$LIN($q$) denote the
problem of satisfying $k$-variable linear equations over ${\bf
Z}_q$, it is known \cite{Hastad:01} that MAX-$k$LIN($q$) has
$(1-\eps, 1/q + \eps)$-hardness for all $k \geq 3$, and this is
sharp. However it seems that current PCP theorems are not strong
enough to give sharp hardness of approximation results for 2-CSPs
(e.g., constraint satisfaction problems on graphs).  The
influential paper of Khot~\cite{Khot:02} introduced the ``Unique
Games Conjecture'' (UGC) in order to make progress on 2-CSPs; UGC
states that a certain 2-CSP over a large domain has $(1-\eps,
\eps)$-hardness.

Interestingly, it seems that using UGC to prove hardness results
for other 2-CSPs typically crucially requires strong results about
influences and noise stability of boolean functions. For example,
\cite{Khot:02}'s analysis of MAX-$2$LIN($2$) required an upper
bound on $\Stab_{1-\eps}(f)$ for small $\eps$ among balanced
functions $f : \bits^n \to \bits$ with small influences; to get
this, Khot used the following deep result of
Bourgain~\cite{Bourgain:02} from 2001:
\begin{theorem}[Bourgain~\cite{Bourgain:02}] \label{thm:bourgain}
If $f : \bits^n \to \bits$ satisfies $\E[f] = 0$ and $\Inf_i(f)
\leq 10^{-d}$ for all $i \in [n]$, then
\[
 \sum_{|S| > d} \hat{f}(S)^2 \geq d^{-1/2 - O(\sqrt{\log \log d / \log d})} = d^{-1/2 -o(1)}.
\]
\end{theorem}
Note that Bourgain's theorem has the following easy corollary:
\begin{corollary} \label{cor:eps1/2-} If $f : \bits^n \to \bits$ satisfies $\E[f] = 0$ and $\Inf_i(f)
\leq 2^{-O(1/\eps)}$ for all $i \in [n]$, then
\[
\Stab_{1-\eps}(f) \leq 1 - \eps^{1/2 + o(1)}.
\]
\end{corollary}
This corollary enabled Khot to show $(1-\eps, 1-\eps^{1/2 +
o(1)})$-hardness for MAX-$2$LIN($2$), which is close to sharp (the
algorithm of Goemans-Williamson~\cite{GoemansWilliamson:95}
achieves $1-O(\sqrt{\eps})$).  As an aside, we note that Khot and
Vishnoi~\cite{KhotVishna:u} recently used
Corollary~\ref{cor:eps1/2-} to prove that negative type metrics do
not embed into $\ell_1$ with constant distortion.

Another example of this comes from the work of~\cite{KKMO:04}.
Among other things,~\cite{KKMO:04} studied the MAX-CUT problem: Given
an undirected graph, partition the vertices into two parts so as
to maximize the number of edges with endpoints in different parts.
The paper introduced the Majority Is Stablest
Conjecture~\ref{conj:MIST} and showed that together with UGC it
implied $(\half + \half \rho - \eps, \half + {\textstyle
\frac{1}{\pi}} \arcsin \rho + \eps)$-hardness for MAX-CUT.  In
particular, optimizing over $\rho$ (taking $\rho \approx .69$) 
implies MAX-CUT is
$.878$-hard to approximate, matching the groundbreaking algorithm
of Goemans and Williamson~\cite{GoemansWilliamson:95}.
%The
%paper~\cite{KKMO:04} also described a ``Plurality Is Stablest
%Conjecture'' for functions $f : [q]^n \to [q]$ and showed this
%would have strong applications for approximate graph coloring and
%UGC itself.

\subsubsection{Consequences of confirming the conjecture}
In Theorem~\ref{thm:MIST} we confirm a generalization of the
Majority Is Stablest conjecture.
%also giving Plurality Is Stablest
We give a slightly simplified statement of this theorem here:
\paragraph{Theorem~\ref{thm:MIST}}\emph{Let $f : \Omega_1 \times \cdots \times \Omega_n \to [0,1]$ be a
function on a discrete product probability space and assume that
for each $i$ the minimum probability of any atom in $\Omega_i$ is
at least $\alpha \leq 1/2$.  Further assume that $\Inf_i(f) \leq
\tau$ for all $i$. Let $\mu = \E[f]$. Then for any $0 \leq \rho <
1$,
\[
\Stab_\rho(f) \leq \lim_{n\to \infty} \Stab_\rho(\Thr^{(\mu)}_n) +
O\Bigl({\textstyle \frac{\log \log (1/\tau)}{\log(1/\tau)}}\Bigr),
\]
where $\Thr^{(\mu)}_n : \bits^n \to \{0,1\}$ denotes the symmetric
threshold function with expectation closest to $\mu$, and the
$O(\cdot)$ hides a constant depending only on $\alpha$ and $1-\rho$.}\\

%\noindent (We remark that Theorem~\ref{thm:MIST} in fact only
%needs $f$ to have small ``low-degree influences'' as defined
%Definition~\ref{def:low-degree-influence}, a distinction crucial
%for PCP applications.)\\

We now give some consequences of this theorem:
\begin{theorem} In the terminology of Kalai~\cite{Kalai:02}, any odd,
balanced social choice function $f$ with either
\begin{itemize}
\item $o_n(1)$ influences or
\item such that $f$ is transitive
\end{itemize}
has probability at most $3/4 + (3/2\pi)\arcsin(1/3) + o_n(1) \approx
.9123$ of producing a rational outcome. The majority function on
$n$ inputs achieves this bound, $3/4 + (3/2\pi)\arcsin(1/3) +
o_n(1)$.
\end{theorem}

By looking at the series expansion of $\frac{2}{\pi} \arcsin(1-\eps)$
we obtain the following strengthening of Corollary~\ref{cor:eps1/2-}.
\begin{corollary} \label{cor:eps1/2} If $f : \bits^n \to \bits$ satisfies $\E[f] = 0$ and $\Inf_i(f)
\leq \eps^{-O(1/\eps)}$ for all $i \in [n]$, then
\[
\Stab_{1-\eps}(f) \leq 1 - ({\textstyle \frac{\sqrt{8}}{\pi} -
o(1)})\eps^{1/2}.
\]
\end{corollary}
Using Corollary~\ref{cor:eps1/2} instead of
Corollary~\ref{cor:eps1/2-} in Khot~\cite{Khot:02} we obtain
\begin{corollary} MAX-$2$LIN($2$) and MAX-2SAT have $(1-\eps, 1 -
O(\eps^{1/2}))$-hardness. \rnote{Actually, for MAX-2LIN(2) we
probably exactly match (to the constant factor) the algorithm of
GW.  Check?}
\end{corollary}
More generally,~\cite{KKMO:04} now implies
\begin{corollary} MAX-CUT has $(\half + \half \rho - \eps, \half + {\textstyle
\frac{1}{\pi}} \arcsin \rho + \eps)$-hardness for each $\rho$ and
all $\eps > 0$, assuming UGC only.  In particular, the
Goemans-Williamson .878-approximation algorithm is best possible,
assuming UGC only.
\end{corollary}
The following two results are consequences of a generalization of
``Majority is Stablest'' as shown in~\cite{KKMO:04}:
\begin{theorem} UGC implies that for each $\eps > 0$ there exists
$q = q(\eps)$ such that MAX-$2$LIN($q$) has $(1-\eps,
\eps)$-hardness.  Indeed, this statement is \emph{equivalent} to
UGC.
\end{theorem}
\begin{theorem} The MAX-$q$-CUT problem, i.e.~Approximate
$q$-Coloring, has $(1 - 1/q + q^{2+o(1)})$-hardness factor, 
assuming UGC only.
This asymptotically matches the approximation factor obtained by
Frieze and Jerrum~\cite{FriezeJerrum:95}.
\end{theorem}

\subsection{It Ain't Over Till It's Over}
The It Ain't Over Till It's Over conjecture was originally made by
Kalai and Friedgut~\cite{Kalai:01} in 
studying social indeterminacy~\cite{FrKaNa:02,Kalai:04}. 
The setting here is similar to the
setting of Arrow's Theorem from Section~\ref{sec:misc-discuss}
except that there are an arbitrary finite number of candidates.
Let $R$ denote the (asymmetric) relation given on the candidates
when the \emph{monotone} social choice function $f$ is used. Kalai
showed that if $f$ has small influences, then the It Ain't Over
Till It's Over Conjecture implies that \emph{every} possible
relation $R$ is achieved with probability bounded away from $0$.
Since its introduction in 2001, the It Ain't Over Till It's Over
problem has circulated widely in the community studying harmonic
analysis of boolean functions.  The conjecture was given as one of
the top unsolved problems in the field at a workshop at Yale in
late 2004.\\

In Theorem~\ref{thm:aint} we confirm the It Ain't Over Till It's
Over conjecture and generalize it to functions on arbitrary
finite product probability spaces with means bounded away from 0
and 1.  Further, the asymptotics we give show that symmetric
threshold functions (e.g., Majority in the case of mean $1/2$) are
the ``worst'' examples.  We give a slightly simplified statement
of Theorem~\ref{thm:aint} here:
\paragraph{Theorem~\ref{thm:aint}}  \emph{Let $0 < \rho < 1$ and
let $f : \Omega_1 \times \cdots \times \Omega_n \to [0,1]$ be a
function on a discrete product probability space; assume that for
each $i$ the minimum probability of any atom in $\Omega_i$ is at
least $\alpha \leq 1/2$.  Then there exists $\eps(\rho,\mu) > 0$ such that 
if $\eps < \eps(\rho,\mu)$ and $\Inf_i(f) \leq \eps^{O(\sqrt{\log(1/\eps)})}$ for all $i$ and $\mu = \E[f]$ then 
\[
\Pr[V_\rho f > 1 - \delta] \leq \eps
\]
and
\[
\Pr[V_\rho f < \delta] \leq \eps
\]
provided
\[
\delta < \eps^{\rho/(1-\rho) + O(1/\sqrt{\log(1/\eps)})},
\]
where the $O(\cdot)$ hides a constant depending only on $\alpha$,
$1-\mu$, $\rho$, and $1-\rho$.}\\

%\noindent (As with Theorem~\ref{thm:MIST}, Theorem~\ref{thm:aint}
%only needs $f$ to have small ``low-degree influences''.  To derive
%the two-sided bound of Conjecture~\ref{conj:ain't} from
%Theorem~\ref{thm:aint}, simply use the theorem twice --- once with
%$f$, once with $1 - f$.)

\section{The invariance principle} \label{sec:invariance}

\subsection{Setup and notation} \label{sec:setup}

In this section we will describe the setup and notation necessary
for our invariance principle.  Recall that we are interested in
functions on finite product probability spaces, $f : \Omega_1
\times \cdots \times \Omega_n \to \R$.  For each $i$, the space of
all functions $\Omega_i \to \R$ can be expressed as the span of a
finite set of orthonormal random variables, $\boldX_{i,0} = 1,
\boldX_{i,1}, \boldX_{i, 2}, \boldX_{i,3}, \dots$; then $f$ can be
written as a multilinear polynomial in the $\boldX_{i,j}$'s. In
fact, it will be convenient for us to mostly disregard the
$\Omega_i$'s and work directly with sets of orthonormal random
variables; in this case, we can even drop the restriction of
finiteness.  We thus begin with the following definition:

\begin{definition}
We call a collection of finitely many orthonormal real random
variables, one of which is the constant $1$, an \emph{orthonormal
ensemble}.  We will write a typical \emph{sequence} of $n$
orthonormal ensembles as $\CalX = (\CalX_1, \dots, \CalX_n)$,
where $\CalX_i = \{\boldX_{i,0} = 1, \boldX_{i,1}, \dots,
\boldX_{i,m_i}\}$.  We call a sequence of orthonormal ensembles
$\CalX$ \emph{independent} if the ensembles are independent families
of random variables.

We will henceforth be concerned only with independent sequences of
orthonormal ensembles, and we will call these \emph{sequences of
ensembles}, for brevity.
\end{definition}

\begin{remark}  \label{rem:simple} Given a sequence of independent random variables
$\boldX_1, \dots, \boldX_n$ with $\E[\boldX_i] = 0$ and
$\E[\boldX_i^2] = 1$ (as in Theorem~\ref{thm:simple}), we can view
them as a sequence of ensembles $\CalX$ by renaming $\boldX_i =
\boldX_{i,1}$ and setting $\boldX_{i,0} = 1$ as required.
\end{remark}

\begin{definition}  We denote by $\CalG$ the \emph{Gaussian sequence of
ensembles}, in which $\CalG_i = \{\boldG_{i,0} = 1, \boldG_{i,1},
\boldG_{i,2}, \dots\}$ and all $\boldG_{i,j}$'s with $j \geq 1$
are independent standard Gaussians.
\end{definition}

As mentioned, we will be interested in \emph{multilinear
polynomials} over sequences of ensembles.  By this we mean sums of
products of the random variables, where each product is obtained
by multiplying one random variable from each ensemble.
\begin{definition}
A \emph{multi-index} $\boldsigma$ is a sequence $(\sigma_1, \dots,
\sigma_n)$ in $\N^n$; the \emph{degree} of $\boldsigma$, denoted
$|\boldsigma|$, is $|\{i \in [n] : \sigma_i > 0\}|$.  Given a
doubly-indexed set of indeterminates $\{x_{i,j}\}_{i \in [n], j
\in \N}$, we write $x_\boldsigma$ for the monomial $\prod_{i =
1}^n x_{i,\sigma_i}$.  We now define a \emph{multilinear
polynomial} over such a set of indeterminates to be any expression
 \begin{equation} \label{eqn:Q}
Q(x) = \sum_{\boldsigma} c_\boldsigma x_\boldsigma
 \end{equation}
where the $c_\boldsigma$'s are real constants, all but finitely
many of which are zero. The \emph{degree} of $Q(x)$ is
$\max\{|\boldsigma| : c_\boldsigma \neq 0\}$, at most $n$.  We
also use the notation
\[
Q^{\leq d}(x) = \sum_{|\boldsigma| \leq d} c_\boldsigma
x_\boldsigma
\]
and the analogous $Q^{= d}(x)$ and $Q^{> d}(x)$.
\end{definition}

Naturally, we will consider applying multilinear polynomials $Q$
to sequences of ensembles $\CalX$; the distribution of these
random variables $Q(\CalX)$ is the subject of our invariance
principle.  Since $Q(\CalX)$ can be thought of as a function on a
product space $\Omega_1 \times \cdots \times \Omega_n$ as
described at the beginning of this section, there is a consistent
way to define the notions of influences, $T_\rho$, and noise
stability from Section~\ref{sec:general}.  For example, the
``influence of the $i$th ensemble on $Q$'' is
\[
\Inf_i(Q(\CalX)) = \E[\Var[Q(\CalX) \mid \CalX_1, \dots,
\CalX_{i-1}, \CalX_{i+1}, \dots, \CalX_n]].
\]
Using independence and orthonormality, it is easy to show the
following formulas, familiar from harmonic analysis of boolean
functions:
\begin{proposition} \label{prop:infQ}
Let $\CalX$ be a sequence of ensembles and $Q$ a multilinear
polynomial as in~(\ref{eqn:Q}). Then
\[
 \E[Q(\CalX)] = c_{\bf{0}}; \qquad \E[Q(\CalX)^2] = \sum_{\boldsigma}
 c_\boldsigma^2; \qquad
 \Var[Q(\CalX)] = \sum_{|\boldsigma| > 0} c_\boldsigma^2;
\]
\[
 \Inf_i(Q(\CalX)) = \sum_{\boldsigma : \sigma_i > 0} c_\boldsigma^2;
 \qquad T_\rho Q(\CalX) = \sum_{\boldsigma} \rho^{|\boldsigma|} c_\boldsigma
 \CalX_\boldsigma; \qquad \Stab_\rho(Q(\CalX)) = \sum_{\boldsigma} \rho^{|\boldsigma|}
 c_\boldsigma^2.
\]
\end{proposition}
Note that in each case above, the formula does not depend on the
sequence of ensembles $\CalX$; it only depends on $Q$.  Thus we
are justified in henceforth writing $\E[Q]$, $\E[Q^2]$, $\Var[Q]$,
$\Inf_i(Q)$, and $\Stab_\rho(Q)$, and in treating $T_\rho$ as a
formal operator on multilinear polynomials:
\begin{definition} \label{def:T_poly}
For $\rho \in [0,1]$ we define the operator $T_\rho$ as acting
formally on multilinear polynomials $Q(x)$ as in~(\ref{eqn:Q}) by
\[
(T_\eta Q)(x) = \sum_{\boldsigma} \rho^{|\boldsigma|} c_\boldsigma
x_\boldsigma.
\]
\end{definition}
Note that for every sequence of ensembles, we have that
Definition~\ref{def:T_poly} agrees with Definition~\ref{def:T_general}.

\bigskip
We end this section with a short discussion of ``low-degree
influences'', a notion that has proven crucial in the analysis of
PCPs (see, e.g., \cite{KKMO:04}).
%\begin{definition}
%Given a multilinear polynomial $Q(x)$ as in~(\ref{eqn:Q}), we define
%the multilinear polynomial
%\[
%Q^{\leq d}(x) = \sum_{|\boldsigma| \leq d} c_\boldsigma
%x_\boldsigma
%\]
%and analogously $Q^{= d}(x)$ and $Q^{> d}(x)$.
%\end{definition}
\begin{definition} \label{def:low-degree-influence} The \emph{$d$-low-degree influence of the $i$th ensemble on $Q(\CalX)$}
is
\[
\Inf^{\leq d}_i(Q(\CalX)) = \Inf^{\leq d}_i(Q) = \sum_{\boldsigma:
|\boldsigma| \leq d, \sigma_i > 0} c_\boldsigma^2.
\]
Note that this gives a way to define low-degree influences
$\Inf^{\leq d}_i(f)$ for functions $f : \Omega_1 \times \cdots
\Omega_n \to \R$ on finite product spaces.
\end{definition}
There isn't an especially natural interpretation of $\Inf_i^{\leq
d}(f)$.  However, the notion is important for PCPs due to the fact
that a function with variance $1$ cannot have too many coordinates
with substantial low-degree influence; this is reflected in the
following easy proposition:
\begin{proposition} \label{prop:infD}
Suppose $Q$ is multilinear polynomial as in~(\ref{eqn:Q}).  Then
\[
\sum_i \Inf_i^{\leq d}(Q) \leq d \cdot \Var[Q].
\]
\end{proposition}

\ignore{
\begin{proposition} $\Stab_\rho(Q) = \E[(T_{\sqrt{\rho}} Q)^2]$.
\end{proposition}
} \ignore{
\begin{remark}
Note that for a general sequence of ensembles $\CalX$ it holds that
$(T_\eta Q)(x)$ is the expected value of $Q(y)$, where for each
ensemble $i$ independently with probability $\eta$ it holds that
$x_{i,j} = y_{i,j}$ for all $j$ and with probability $1-\eta$ the
vector $(y_{i,j})_j$ is drawn from $\CalX_i$.

Note furthermore, that for $\CalG$, the operator $T_{\eta}$ is the
usual Ornstein-Uhlenbeck operator
\[
(T_\rho f)(x) = \Ex_{y}[f(\rho x + \sqrt{1-\rho^2}\, y)],
\]
where the expected value is with respect to the Gaussian measure.
\end{remark}
}

\subsection{Hypercontractivity}
As mentioned in Section~\ref{sec:intro-invariance}, our invariance
principle requires that the ensembles involved to be
hypercontractive in a certain sense.  Recall that a random
variable $\boldY$ is said to be ``$(p,q,\eta)$-hypercontractive''
for $1 \leq p \leq q < \infty$ and $0 < \eta < 1$ if
\begin{equation} \label{eqn:easy-hc}
\|a + \eta \boldY\|_q \leq \|a + \boldY\|_p
\end{equation}
for all $a \in \R$.  This type of hypercontractivity was
introduced (with slightly different notation)
in~\cite{KrakowiakSzulga:88}. Some basic facts about
hypercontractivity are explained in Appendix~\ref{app:hc}; much
more can be found in \cite{KwapienWoyczynski:92}. Here we just
note that for $q>2$ a random variable $\boldY$ is
$(2,q,\eta)$-hypercontractive with some $\eta \in (0,1)$ if and
only if $\E[\boldY]=0$ and $\E[|\boldY|^{q}]<\infty.$ Also, if
$\boldY$ is $(2,q,\eta)$-hypercontractive then $\eta \leq
(q-1)^{-1/2}.$\\

We now define our extension of the notion of hypercontractivity to
sequences of ensembles:

\begin{definition}
Let $\CalX$ be a sequence of ensembles.  For $1 \leq p \leq q <
\infty$ and $0 < \eta < 1$ we say that $\CalX$ is
\emph{$(p,q,\eta)$-hypercontractive} if
\[
\|(T_\eta Q)(\CalX)\|_q \leq \|Q(\CalX)\|_p
\]
for every multilinear polynomial $Q$ over $\CalX$.
\end{definition}
Since $T_{\eta}$ is a contractive semi-group, we have
\begin{remark} If $\CalX$ is $(p,q,\eta)$-hypercontractive then it
is $(p,q,\eta')$-hypercontractive for any $0 < \eta' \leq \eta$.
\end{remark}

There is a related notion of hypercontractivity for \emph{sets} of
random variables which considers all polynomials in the variables,
not just multilinear polynomials; see, e.g.,
Janson~\cite{Janson:97}.  Several of the properties of this notion
of hypercontractivity carry over to our setting of sequences of
ensembles.  In particular, the following facts can easily be
proved by repeating the analogous proofs in~\cite{Janson:97}; for
completeness, we give the proofs in Appendix~\ref{app:hc}.

\begin{proposition} \label{prop:join-hypercon}
Suppose $\CalX$ is a sequence of $n_1$ ensembles and $\CalY$ is an
independent sequence of $n_2$ ensembles.  Assume both are
$(p,q,\eta)$-hypercontractive.  Then the sequence of ensembles
$\CalX \cup \CalY = (\CalX_1, \dots, \CalX_{n_1}, \CalY_1, \dots,
\CalY_{n_2})$ is also $(p,q,\eta)$-hypercontractive.
\end{proposition}

\begin{proposition} \label{prop:hypercon}
Let $\CalX$ be a $(2,q,\eta)$-hypercontractive sequence of
ensembles and $Q$ a multilinear polynomial over $\CalX$ of degree
$d$.  Then
\[
\|Q(\CalX)\|_q \leq \eta^{-d} \; \|Q(\CalX)\|_2.
\]
\end{proposition}

In light of Proposition~\ref{prop:join-hypercon}, to check that a
sequence of ensembles is $(p,q,\eta)$-hypercontractive it is
enough to check that each ensemble individually is
$(p,q,\eta)$-hypercontractive (as a ``sequence'' of length 1); in
turn, it is easy to see that this is equivalent to checking that
for each $i$, all linear combinations of the random variables
$\boldX_{i,1}, \dots, \boldX_{i, m_i}$ are hypercontractive in the
traditional sense of~(\ref{eqn:easy-hc}).\\

We end this section by recording the optimal hypercontractivity
constants for the ensembles we consider.  The result for $\pm 1$
Rademacher variables is well known and due originally to
Bonami~\cite{Bonami:70} and independently
Beckner~\cite{Beckner:75}; the same result for Gaussian and uniform
random variables is also well known and in fact follows easily from the
Rademacher case. The optimal hypercontractivity constants for general
finite spaces was recently determined by Wolff~\cite{Wolff:u}
(see also~\cite{Oleszkiewicz:03}):
\begin{theorem} \label{thm:bonami} Let $\boldX$ denote either a uniformly random $\pm 1$
bit, a standard one-dimensional Gaussian, or a random variable
uniform on $[-\sqrt{3}, \sqrt{3}]$. Then $\boldX$ is $(2, q,
(q-1)^{-1/2})$-hypercontractive.
\end{theorem}
\begin{theorem} \label{thm:wolff} (Wolff)\ \ Let $\boldX$ be any
mean-zero random variable on a finite probability space in which
the minimum nonzero probability of any atom is $\alpha \leq 1/2$.
Then $\boldX$ is $(2, q, \eta_q(\alpha))$-hypercontractive, where
\[
\eta_q(\alpha) = \left(\frac{A^{1/q'} - A^{-1/q'}}{A^{1/q} -
A^{-1/q}} \right)^{-1/2}
\]
\[
\text{with } \quad A = \frac{1-\alpha}{\alpha}, \quad 1/q + 1/q' =
1.
\]
\end{theorem}
Note the following special case:
\begin{proposition} \label{prop:wolff}
\[
\eta_3(\alpha) = \left(A^{1/3} + A^{-1/3}\right)^{-1/2} \;\;
\mathop{\sim}^{\alpha \to 0} \quad \alpha^{1/6},
\]
and also
\[
\half \alpha^{1/6} \leq \eta_3(\alpha) \leq 2^{-1/2},
\]
for all $\alpha \in [0,1/2]$.
\end{proposition}
For general random variables with bounded moments we have the
following results, proved in Appendix~\ref{app:hc}:
\begin{proposition} \label{prop:bdd}
Let $\boldX$ be a mean-zero random variable satisfying
$\E[|\boldX|^{q}]< \infty$.  Then $\boldX$ is
$(2,q,\eta_{q})$-hypercontractive with $\eta_{q}=\frac{\| \boldX
\|_{2}}{2\sqrt{q-1}\| \boldX \|_{q}}.$
\end{proposition}
In particular, when $\E[\boldX] = 0$, $\E[\boldX^{2}]=1$, and
$\E[|\boldX|^{3}] \leq \beta$, we have that $\boldX$ is
$(2,3,2^{-3/2}\beta^{-1/3})$-hypercontractive.
\begin{proposition} \label{prop:add}
Let $\boldX$ be a mean-zero random variable satisfying
$\E[|\boldX|^{q}]< \infty$ and let $\boldV$ be a random variable
independent of $\boldX$ with $\P[\boldV=0]=1-\rho$ and
$\P[\boldV=1]=\rho.$ Then $\boldV\boldX$ is
$(2,q,\xi_{q})$-hypercontractive with $\xi_{q}=\frac{\| \boldX
\|_{2}}{2\sqrt{q-1}\| \boldX \|_{q}} \cdot
\rho^{\frac{1}{2}-\frac{1}{q}}.$
\end{proposition}

\subsection{Hypotheses for invariance theorems --- some families of ensembles}
All of the variants of our invariance principle that we prove in
this section will have similar hypotheses.  Specifically, they
will be concerned with a multilinear polynomial $Q$ over two
hypercontractive sequences of ensembles, $\CalX$ and $\CalY$;
furthermore, $\CalX$ and $\CalY$ will be assumed to have satisfy a
``matching moments'' condition, as described below.  We will now
lay out four hypotheses --- $\hone,\htwo, \hthree$, and $\hfour$
that will be used in the theorems of this section. As can easily
be seen (using Theorems~\ref{thm:bonami} and~\ref{thm:wolff} and
Proposition~\ref{prop:wolff}; see also Appendix~\ref{app:hc}), the
hypothesis $\hone$ generalizes $\htwo, \hthree$, and $\hfour$;
hence all proofs will be carried out only in the setting of
$\hone$. However the amount of notation and number of parameters
under $\hone$ is quite cumbersome, and the reader who is
interested mainly in functions on finite product spaces
($\hthree$) or just boolean functions where $\bits^n$ has the
uniform distribution ($\hfour$) may find it easier to proceed
through the
proofs and results in the restricted cases.\\
%The hypothesis
%$\hfour$ will be used to derive a multilinear
%version of the Berry-Esseen inequality (Theorem~\ref{thm:simple}).\\

Herewith our hypotheses:

\begin{enumerate}
\item[$\hone$]
Let $r \geq 3$ be an integer and let $\CalX$ and $\CalY$ be
independent sequences of $n$ ensembles which are
$(2,r,\eta)$-hypercontractive; recall that $\eta \leq
(r-1)^{-1/2}$. Assume furthermore that for all $1 \leq i \leq n$
and all sets $\Sigma \subset \N$ with $|\Sigma| < r$, the
sequences $\CalX$ and $\CalY$ satisfy the ``matching moments''
condition
\begin{equation} \label{eq:mixed_r_moments}
 \E \left[\prod_{\sigma \in \Sigma} \boldX_{i,\sigma} \right] =
 \E \left[\prod_{\sigma \in \Sigma} \boldY_{i,\sigma} \right].
\end{equation}
Finally, let $Q$ be a multilinear polynomial as
in~(\ref{eqn:Q}).\\

We remark that in $\hone$, if $r = 3$ then the matching moment
conditions hold automatically since the sequences are orthonormal.
We also remark that we have added the condition $\eta \leq
(r-1)^{-1/2}$ so that we can take $\CalY = \CalG$, the Gaussian
sequence of ensembles (see Theorem~\ref{thm:bonami}).

\item[$\htwo$] Let $r = 3.$ Let $\CalX$ and $\CalY$ be independent
sequences of ensembles in which each ensemble has only two random
variables, $\boldX_{i,0} = 1$ and $\boldX_{i,1} = \boldX_i$
(respectively, $\boldY_{i,0} = 1$, $\boldY_{i,1} = \boldY_i$), as
in Remark~\ref{rem:simple}.  Further assume that each $\boldX_i$
(respectively $\boldY_i$) satisfies $\E[\boldX_i]=0,$
$\E[\boldX_i^{2}]=1$ and $\E[|\boldX_i|^{3}] \leq \beta.$ Put
$\eta=2^{-3/2}\beta^{-1/3}$, so $\CalX$ and $\CalY$ are
$(2,3,\eta)$-hypercontractive. Finally, let $Q$ be a multilinear
polynomial as in~(\ref{eqn:Q}).\\

The hypothesis $\htwo$ is used to derive the multilinear version
of the Berry-Esseen inequality given in Theorem~\ref{thm:simple}.

\item[$\hthree$] Let $r = 3$ and let $\CalX$ be a sequence of $n$
ensembles in which the random variables in each ensemble $\CalX_i$
form a basis for the real-valued functions on some finite
probability space $\Omega_i$.  Further assume that the least
nonzero probability of any atom in any $\Omega_i$ is $\alpha \leq
1/2$, and let $\eta = \half \alpha^{1/6}$.  Let $\CalY$ be any
independent $(2,3,\eta)$-hypercontractive sequence of ensembles.
Finally, let $Q$ be a multilinear polynomial as
in~(\ref{eqn:Q}).\\

We remark that $Q(\CalX)$ in $\hthree$ encompasses \emph{all}
real-valued functions $f$ on finite product spaces, including
the familiar cases of the $p$-biased discrete cube (for which
$\alpha = \min\{p, 1-p\}$) and the set $[q]^n$ with uniform
measure (for which $\alpha = 1/q$).  Note also that $\eta \leq
2^{-1/2}$ so we may take $\CalY$ to be the Gaussian sequence of
ensembles.

\item[$\hfour$] Let $r = 4$ and $\eta = 3^{-1/2}$.  Let $\CalX$
and $\CalY$ be independent sequences of ensembles in which each
ensemble has only two random variables, $\boldX_{i,0} = 1$ and
$\boldX_{i,1} = \boldX_i$ (respectively, $\boldY_{i,0} = 1$,
$\boldY_{i,1} = \boldY_i$), as in Remark~\ref{rem:simple}. Further
assume that each $\boldX_i$ (respectively $\boldY_i$) is either a)
a uniformly random $\pm 1$ bit; b) a standard one-dimensional
Gaussian; or c) uniform on $[-3^{1/2}, 3^{1/2}]$. Hence $\CalX$
and $\CalY$ are $(2,4,\eta)$-hypercontractive. Finally, let $Q$ be
a multilinear polynomial as
in~(\ref{eqn:Q}).\\

Note that this simplest of all hypotheses allows for arbitrary
real-valued functions on the uniform-measure discrete cube $f :
\bits^n \to \R$.  Also, under $\hfour$, $Q$ is just a multilinear
polynomial in the usual sense over the $\boldX_i$'s or
$\boldY_i$'s; in particular, if $f : \bits^n \to \R$ then $Q$ is
the ``Fourier expansion'' of $f$.  Finally, note that the matching
moments condition~(\ref{eq:mixed_r_moments}) holds in $\hfour$
since it requires $\E[X_t^3] = \E[Y_t^3]$ for each $t$, and this
is true since both equal $0$.
\end{enumerate}

\subsection{Basic invariance principle, $\CalC^r$ functional version}
The essence of our invariance principle is that if $Q$ is of
bounded degree and has low influences then the random variables
$Q(\CalX)$ and $Q(\CalY)$ are close in distribution. The simplest
way to formulate this conclusion is to say that if $\TestFcn : \R
\to \R$ is a sufficiently nice ``test function'' then
$\TestFcn(Q(\CalX))$ and $\TestFcn(Q(\CalY))$ are close in
expectation.

\begin{theorem} \label{thm:cj}  Assume hypothesis $\hone,\htwo,\hthree$, or
$\hfour.$ Further assume $\Var[Q] \leq 1$, $\deg(Q) \leq d$, and $\Inf_i(Q)
\leq \tau$ for all $i$.  Let $\TestFcn : \R \to \R$ be a $\CalC^r$
function with $|\TestFcn^{(r)}| \leq B$ uniformly.  Then
\[
\Bigl| \E\bigl[\TestFcn(Q(\CalX))\bigr] -
\E\bigl[\TestFcn(Q(\CalY))\bigr] \Bigr| \leq \eps,
\]
where
\[
\eps = \left\{ \begin{array}{ll}
                   (2B/r!)\,d\,\eta^{-rd}\;\tau^{r/2 - 1} & \text{under } \hone, \\
                   B\,30^d\;\beta^{d}\;\tau^{1/2} & \text{under } \htwo, \\
                   B\,(10\alpha^{-1/2})^d\;\tau^{1/2} & \text{under } \hthree, \\
                   B\,10^d\;\tau & \text{under } \hfour.
               \end{array}   \right.
\]
\end{theorem}

As will be the case in all of our theorems, the results under
$\htwo,\hthree$ and $\hfour$ are immediate corollaries of the
result under $\hone$; one only needs to substitute in $r=3$,
$\eta=2^{-3/2}\beta^{-1/3}$ or $r = 3$, $\eta = \half
\alpha^{1/6}$ or $r = 4$, $\eta = 3^{-1/2}$ (we have also here
used that $(1/3)\,d\,2^{9d/2}$ is at most $30^d$ and that
$(1/3)\,d\,8^d$ and $(1/12)\,d\,9^d$ are at most $10^d$). Thus it
will suffice for us to carry out the proof under $\hone$.\\

%For this
%particular theorem we should point out that the proof becomes
%notationally far simpler under $\hthree$, when the
%ensembles contain only one random variable in addition to 1; yet,
%all of the essential ideas remain present.  Thus we encourage the
%reader to try following the
%proof under \textbf{Hypothesis III} on first reading.\\

\begin{proof}
We begin by defining intermediate sequences between $\CalX$ and
$\CalY$.  For $i = 0, 1, \dots, n$, let $\CalZ^{(i)}$ denote the
sequence of $n$ ensembles $(\CalY_1, \dots, \CalY_i, \CalX_{i+1},
\dots, \CalX_n)$ and let $\boldQ^{(i)} = Q(\CalZ^{(i)})$.  Our
goal will be to show
\begin{equation} \label{eqn:bound}
\Bigl|\E\bigl[\TestFcn(\boldQ^{(i-1)})\bigr] -
\E\bigl[\TestFcn(\boldQ^{(i)})\bigr]\Bigr| \leq
\left(\frac{2B}{r!}\,\eta^{-rd}\right)\cdot\Inf_i(Q)^{r/2}
\end{equation}
for each $i \in [n]$.  Summing this over $i$ will complete the
proof since $\CalZ^{(0)} = \CalX$, $\CalZ^{(n)} = \CalY$, and
\[
\sum_{i=1}^n \Inf_i(Q)^{r/2} \leq \tau^{r/2-1} \cdot \sum_{i=1}^n
\Inf_i(Q) = \tau^{r/2-1} \cdot \sum_{i=1}^n \Inf_i^{\leq d}(Q)
\leq d \tau^{r/2-1},
\]
where we used Proposition~\ref{prop:infD} and $\Var[Q] \leq 1$.\\

\newcommand{\Qt}{\tilde{{\boldQ}}} \newcommand{\Rt}{{\boldR}} \newcommand{\St}{{\boldS}}
Let us fix a particular $i \in [n]$ and proceed to
prove~(\ref{eqn:bound}).  Given a multi-index $\boldsigma$, write
$\boldsigma \setminus i$ for the same multi-index except with
$\sigma_i = 0$.  Now write
\begin{eqnarray*}
\Qt &=& \sum_{\boldsigma : \sigma_i = 0} c_\boldsigma
\CalZ^{(i)}_\boldsigma, \\
\Rt &=& \sum_{\boldsigma : \sigma_i > 0} c_\boldsigma
\boldX_{i,\sigma_i} \cdot \CalZ^{(i)}_{\boldsigma \setminus i}, \\
\St &=& \sum_{\boldsigma : \sigma_i > 0} c_\boldsigma
\boldY_{i,\sigma_i} \cdot \CalZ^{(i)}_{\boldsigma \setminus i}.
\end{eqnarray*}
Note that $\Qt$ and the variables $\CalZ^{(i)}_{\boldsigma
\setminus
  i}$ are independent of the variables in $\CalX_i$ and $\CalY_i$ and
that $\boldQ^{(i-1)} = \Qt + \Rt$ and $\boldQ^{(i)} = \Qt + \St$.
\\

To bound the left side of~(\ref{eqn:bound}) --- i.e.,
$|\E[\TestFcn(\Qt + \Rt) - \TestFcn(\Qt + \St)]|$ --- we use
Taylor's theorem: for all $x, y \in \R$,
\[
\Bigl|\TestFcn(x+y) - \sum_{k=0}^{r-1}
\frac{\TestFcn^{(k}(x)\;y^k}{k!}\Bigr| \leq \frac{B}{r!}\,|y|^r.
\]
In particular,
\begin{equation} \label{eq:R_taylor}
\Bigl|\E[\TestFcn(\Qt + \Rt)] - \sum_{k=0}^{r-1}
\E\Bigl[\frac{\TestFcn^{(k)}(\Qt)\;\Rt^k}{k!}\Bigr]\Bigr| \leq
\frac{B}{r!}\,\E\bigl[|\Rt|^r\bigr]
\end{equation}
and similarly,
\begin{equation} \label{eq:S_taylor}
\Bigl|\E[\TestFcn(\Qt + \St)] - \sum_{k=0}^{r-1}
\E\Bigl[\frac{\TestFcn^{(k)}(\Qt)\;\St^k}{k!}\Bigr]\Bigr| \leq
\frac{B}{r!}\,\E\bigl[|\St|^r\bigr].
\end{equation}
We will see below that that $\Rt$ and $\St$ have finite $r$ moments.
Moreover,
for $0 \leq k
\leq r$ it holds that $|\TestFcn^{(k)}(\Qt)\,\Rt^k| \leq
|k!\,B\,\Qt^{r-k}\,\Rt^k|$ (and similarly for $\St$). Thus all
moments above are finite.
We now claim
that for all $0 \leq k < r$ it holds that
\begin{equation} \label{eq:S_T_moments}
\E[\TestFcn^{(k)}(\Qt)\,\Rt^k] = \E[\TestFcn^{(k)}(\Qt)\,\St^k].
\end{equation}
Indeed,
\begin{eqnarray} \nonumber
\E[\TestFcn^{(k)}(\Qt)\,\Rt^k] &=& \E \Bigl[\TestFcn^{(k)}(\Qt)
\sum_{\substack{(\boldsigma^1,\dots,\boldsigma^k) \\ \text{s.t. }
\forall t,\;\;\sigma^t_i > 0}} \prod_{t=1}^k c_{\boldsigma^t}
\prod_{t=1}^k \boldX_{i,\sigma^t_i}
\prod_{t=1}^k \CalZ^{(i)}_{\boldsigma^t \setminus i} \Bigr] \\
&=& \label{eq:clt_ind}
\sum_{\substack{(\boldsigma^1,\dots,\boldsigma^k) \\ \text{s.t. }
\forall t,\;\;\sigma^t_i > 0}} \prod_{t=1}^k c_{\boldsigma^i}
\cdot \E \Bigl[\TestFcn^{(k)}(\Qt)
    \prod_{t=1}^k \CalZ^{(i)}_{\boldsigma^t \setminus i} \Bigr]
\cdot \E\Bigl[\prod_{t=1}^k \boldX_{i,\sigma^t_i}\Bigr] \\
 &=& \label{eq:clt_same_dist}
\sum_{\substack{(\boldsigma^1,\dots,\boldsigma^k) \\ \text{s.t. }
\forall t,\;\;\sigma^t_i > 0}} \prod_{t=1}^k c_{\boldsigma^i}
\cdot \E \Bigl[\TestFcn^{(k)}(\Qt)
    \prod_{t=1}^k \CalZ^{(i)}_{\boldsigma^t \setminus i} \Bigr]
\cdot \E\Bigl[\prod_{t=1}^k \boldY_{i,\sigma^t_i}\Bigr] \\
&=& \nonumber \E \left[\TestFcn^{(k)}(\Qt)\,\St^k \right].
\end{eqnarray}

\noindent The equality in (\ref{eq:clt_ind}) follows since
$\CalZ^{(i)}_{\boldsigma^t \setminus i}$ and $\Qt$ are independent
of the variables in $\CalX_i$ and $\CalY_i$. The equality in
(\ref{eq:clt_same_dist}) follows from the matching moments
condition~(\ref{eq:mixed_r_moments}).\\

From (\ref{eq:R_taylor}), (\ref{eq:S_taylor}) and
(\ref{eq:S_T_moments}) it follows that
\begin{equation} \label{eq:taylor_bd}
 |\E[\TestFcn(\Qt + \Rt) - \TestFcn(\Qt + \St)]| \leq \frac{B}{r!}\,(\E[|\Rt|^r] +
 \E[|\St|^r]).
\end{equation}
We now use hypercontractivity.  By
Proposition~\ref{prop:join-hypercon} each $\CalZ^{(i)}$ is
$(2,r,\eta)$-hypercontractive.  Thus by
Proposition~\ref{prop:hypercon},
\begin{equation} \label{eq:R_S_hyper}
\E[|\Rt|^r] \leq \eta^{-rd} \E[\Rt^2]^{r/2}, \quad
\E[|\St|^r] \leq \eta^{-rd} \E[\St^2]^{r/2}.
\end{equation}
However,
\begin{equation} \label{eq:S_R_inf}
\E[\St^2] = \E[\Rt^2] = \sum_{\boldsigma : \sigma_i > 0}
c_\boldsigma^2 = \Inf_i(Q).
\end{equation}
Combining (\ref{eq:taylor_bd}), (\ref{eq:R_S_hyper}) and
(\ref{eq:S_R_inf}) it follows that
\[
 |\E[\TestFcn(\Qt + \Rt) - \TestFcn(\Qt + \St)]| \leq
\left(\frac{2B}{r!}\,\eta^{-rd}\right)\cdot\Inf_i(Q)^{r/2}
\]
confirming~(\ref{eqn:bound}) and completing the proof.
\end{proof}

\subsection{Invariance principle --- other functionals, and smoothed version}

Our basic invariance principle shows that $\E[\TestFcn(Q(\CalX))]$
and $\E[\TestFcn(Q(\CalY))]$ are close if $\TestFcn$ is a
$\CalC^r$ functional with bounded $r$th derivative.  To show that
the distributions of $Q(\CalX)$ and $Q(\CalY)$ are close in other
senses we need the invariance principle for less smooth
functionals.  This we can obtain using straightforward
approximation arguments; we defer the proof of Theorem~\ref{thm:supertheorem}
which follows to Section~\ref{sec:proofs}.\\

Theorem~\ref{thm:supertheorem} shows closeness of distribution in
two senses.  The first is closeness in \emph{L\'{e}vy's metric};
recall that the distance between two random variables $\boldR$ and
$\boldS$ in L\'{e}vy's metric is
 \[
 d_L(\boldR, \boldS) = \inf\{\lam > 0 : \quad
\forall t \in \R,\;\; \Pr[\boldS \leq t - \lam] - \lam \leq
\Pr[\boldR \leq t] \leq \Pr[\boldS \leq t + \lam] + \lam\}.
 \]
We also show the distributions are close in the usual sense with a
weaker bound; the proof of this goes by comparing the
distributions of $Q(\CalX)$ and $Q(\CalY)$ to $Q(\CalG)$ and
noting that bounded-degree Gaussian polynomials are known to have
low ``small ball probabilities''.  Finally,
Theorem~\ref{thm:supertheorem} also shows $L^1$ closeness and, as
a technical necessity for applications, shows closeness under the
functional $\trunc : \R \to \R$ defined by
\begin{equation} \label{eq:def_trunc}
\trunc(x) = \left\{ \begin{array}{ll}
                    x^2 & \mbox{ if } x \leq 0,\\
                    0 & \mbox{ if } x \in [0,1], \\
                    (x-1)^2 & \mbox{ if } x \geq 1;
                    \end{array} \right.
\end{equation}
this functional gives the squared distance to the interval
$[0,1]$.\\

\begin{theorem} \label{thm:supertheorem}  Assume Hypothesis
$\hone,\htwo,\hthree$, or $\hfour$. Further assume $\Var[Q] \leq
1$, $\deg(Q) \leq d$ and $\Inf_i(Q) \leq \tau$ for all $i$.  Then
\begin{eqnarray}
 \Bigl| \|Q(\CalX)\|_1 - \|Q(\CalY)\|_1 \Bigr| &\leq& O(\eps^{1/r}), \label{eq:superl1} \\
 d_L(Q(\CalX), Q(\CalY)) &\leq& O(\eps^{1/(r+1)}), \label{eq:superlevi} \\
 \Bigl| \E\bigl[\trunc(Q(\CalX))\bigr] - \E\bigl[\trunc(Q(\CalY))\bigr] \Bigl| & \leq &
 O(\eps^{2/r}),
 \label{eq:superl2}
\end{eqnarray}
where $O(\cdot)$ hides a constant depending only on $r$, and
\[
\eps = \left\{ \begin{array}{ll}
                   d\,\eta^{-rd}\;\tau^{r/2 - 1} & \text{under } \hone,\\
                   30^{d}\beta^{d}\;\tau^{1/2} & \text{under } \htwo,\\
                   (10\alpha^{-1/2})^d\;\tau^{1/2} & \text{under } \hthree, \\
                   10^d\;\tau & \text{under } \hfour.\\
               \end{array}   \right.
\]
If in addition $\Var[Q] = 1$ then
 \begin{equation} \label{eq:lim_dist}
   \sup_t\;\Bigl|\P\bigl[Q(\CalX) \leq t\bigr] - \P\bigl[Q(\CalY) \leq t\bigr]\Bigr| \leq
   O\bigl(d\,\eps^{1/(rd+1)}\bigr).
 \end{equation}
\end{theorem}

\bigskip

As discussed in Section~\ref{sec:intro-invariance},
Theorem~\ref{thm:supertheorem} has the unavoidable deficiency of
having error bounds depending on the degree $d$ of $Q$.  This can
be overcome if we first ``smooth'' $Q$ by applying $T_{1 -
\gamma}$ to it, for some $0 < \gamma < 1$.
Theorem~\ref{thm:smooththeorem} which follows will be our main
tool for applications; its proof is a straightforward degree
truncation argument which we also defer to
Section~\ref{sec:proofs}.  As an additional benefit of this
argument, we will show that $Q$ need only have small
\emph{low-degree influences}, $\Inf_i^{\leq d}(Q)$, as opposed to
small influences.  As discussed at the end of
Section~\ref{sec:setup}, this feature has proven essential for
applications involving PCPs.

\begin{theorem} \label{thm:smooththeorem}
Assume hypothesis $\hone, \hthree$,or $\hfour$. Further assume
$\Var[Q] \leq 1$ and \linebreak $\Inf_i^{\leq\,\log(1/\tau)/K}(Q)
\leq \tau \leq$ for all $i$, where
\[
 K = \left\{ \begin{array}{ll}
                   \log(1/\eta) & \text{under } \hone, \\
                   \log(1/\alpha) & \text{under } \hthree, \\
                   1 & \text{under } \hfour.
              \end{array} \right.
\]
Given $0 < \gamma < 1$, write $\boldR = (T_{1-\gamma} Q)(\CalX)$
and $\boldS = (T_{1-\gamma} Q)(\CalY)$. Then
 \begin{eqnarray*}
 d_L(\boldR, \boldS) & \leq &  \tau^{\Omega(\gamma/K)}, \\
 \Bigl|\E\bigl[\trunc(\boldR)\bigr] -
 \E\bigl[\trunc(\boldS)\bigr]\Bigr| &\leq&
 \tau^{\Omega(\gamma/K)},
 \end{eqnarray*}
where the $\Omega(\cdot)$ hides a constant depending only on $r$.

More generally the statement of the theorem holds for 
$\boldR = Q(\CalX), \boldS=Q(\CalY)$ if 
$\Var[Q^{> d}] \leq (1-\gamma)^{2d}$ for all $d$. 
\end{theorem}

\ignore{ Here are the painstaking bounds we can actually get:

 \begin{eqnarray}
   d_L(\boldR, \boldS) &\leq& C_r\,\tau^{1/c_1}\cdot[\log(1/\tau)/\log(1/\rho\eta)]^{1/(r+1)},  \label{eq:levi_smooth} \\
   |\E[\trunc(\boldR)] - \E[\trunc(\boldS)]| &
   \leq & C_r\,\tau^{1/c_2}\cdot[\log(1/\tau)/\log(1/\rho\eta)]^{2/r}, \label{eq:l2_smooth}
 \end{eqnarray}
where
 \begin{equation} \label{eq:c1_c2}
 c_1 = \frac{3r}{r-2} \cdot \frac{\ln(1/\eta)}{\ln(1/\rho)} +
 \frac{2r+2}{r-2},
 \qquad c_2 = \frac{2r}{r-2} \cdot \frac{\ln(1/\eta)}{\ln(1/\rho)} +
 \frac{r}{r-2}.
 \end{equation}
}

\subsection{Proofs of extensions of the invariance principle} \label{sec:proofs}

In this section we will prove Theorems~\ref{thm:supertheorem}
and~\ref{thm:smooththeorem} under hypothesis $\hone$. The results
under $\htwo,\hthree$, and $\hfour$ are corollaries.

\subsubsection{Invariance principle for some $C^0$ and $C^1$
functionals} \label{sec:convolve} In this section we
prove~(\ref{eq:superl1}), (\ref{eq:superlevi}), (\ref{eq:superl2})
of Theorem~\ref{thm:supertheorem}. We do it by approximating the
following functions in the sup norm by smooth functions:
% \begin{eqnarray*}
\[\begin{array}{lll}
     \ell_1(x) = |x|; &
     \Delta_{s,t}(x) = \left\{ \begin{array}{ll}
                    1 & \mbox{ if } x \leq t-s, \\
                    \frac{t-x+s}{2s} & \mbox{ if } x \in [t-s,t+s], \\
                    0 & \mbox{ if } x \geq t+s;
                    \end{array} \right. &
     \trunc(x) = \left\{ \begin{array}{ll}
                    x^2 & \mbox{ if } x \leq 0,\\
                    0 & \mbox{ if } x \in [0,1], \\
                    (x-1)^2 & \mbox{ if } x \geq 1.
                    \end{array} \right.
                  \end{array}
\]

\begin{lemma} \label{lem:approx-functional}
Let $r \geq 2$ be an integer. Then there exist constant $B_r$ for
which the following holds. For all $0 < \lam \leq 1/2$ there exist
$\CalC^\infty$ functions $\ell_1^\lam$, $\Delta_{\lambda,t}^\lam$
and $\trunc^\lam$ satisfying the following:
\begin{itemize}
 \item $\|\ell_1^\lam - \ell_1\|_\infty \leq 2\lam$; and, $\|(\ell_1^\lam)^{(r)}\|_\infty \leq
 4 B_r\,\lambda^{1-r}$.
 \item $\Delta_{\lam,t}^\lam$ agrees with $\Delta_{\lam,t}$ outside
 the interval $(t-2\lam, t+2\lam)$, and is otherwise in $[0,1]$; and, $\|(\Delta_{\lam,t}^\lam)^{(r)}\|_\infty \leq
 B_r\,\lambda^{-r}$.
 \item $\|\trunc^\lam - \trunc\|_\infty \leq 2\lam^2$; and, $\|(\trunc^\lam)^{(r)}\|_\infty \leq
2 B_{r-1} \,\lam^{2-r}$.
\end{itemize}
\end{lemma}

\begin{proof}  Let $f(x) = x 1_{\{x \geq 0\}}$.
We will show that for all $\lam > 0$ there is a $\CalC^\infty$
function $f_\lam$ satisfying the following:
 \begin{itemize}
 \item $f_\lam$ and $f$ agree on $(-\infty, -\lam]$ and $[\lam,
 \infty)$;
 \item $0 \leq f_\lam(x) \leq f(x) + \lam$ on $(-\lam, \lam)$;
 and,
 \item $\|f_\lam^{(r)}\|_\infty \leq 2 B_r\,\lam^{1-r}$.
 \end{itemize}
The construction of $f$ easily gives the construction of the other
functionals by letting $\ell_1^\lambda(x)  =  f_\lam(x) +
f_\lam(-x)$ and
\begin{equation} \label{eq:delta_lam}
 \Delta_{\lambda,t}^\lambda(x)  =  \left\{ \begin{array}{ll}
                \frac{1}{2\lambda}f_\lam(t-x+\lambda) & \mbox{ if } x \geq t,\\
                1 - \frac{1}{2\lambda}f_{\lam}(x-t+\lambda) & \mbox{ if } x \leq
                t; \end{array} \right. \qquad \qquad
 \trunc^\lambda(x) = \left\{ \begin{array}{ll}
                \int_{-\infty}^{x-1} f_\lam(t)dt & \mbox{ if } x \geq 1/2,\\
                \int_{-\infty}^{1-x} f_\lam(t)dt& \mbox{ if } x \leq 0.
                \end{array} \right.
\end{equation}
To construct $f$, first let $\psi$ be a nonnegative $\CalC^\infty$
function satisfying the following: $\psi$ is $0$ outside $(-1,1)$,
$\int_{-1}^1 \psi(x)\,dx = 1$, and $\int_{-1}^1 x \psi(x)\,dx =
0$.
%\begin{itemize}
% \item $\psi$ is $0$ outside $(-1,1)$;
% \item $\int_{-1}^1 \psi(x)\,dx = 1$;
% \item $\int_{-1}^1 x \psi(x)\,dx = 0$.
%\end{itemize}
It is well known that such functions $\psi$ exist.  Define the
constant $B_r$ to be $\|\psi^{(r)}\|_\infty$.\\

Next, write $\psi_\lambda(x) = \psi(x/\lam)/\lam$, so
$\psi_\lambda$ satisfies the same three properties as $\psi$ with
respect to the interval $(-\lambda, \lambda)$ rather than
$(-1,1)$.  Note that $\|\psi_\lambda^{(r)}\|_\infty =
B_r\,\lambda^{-1-r}$.\\

Finally, take $f_\lambda = f * \psi_\lambda$, which is
$\CalC^\infty$. The first two properties demanded of $f$ follow
easily.  To see the third, first note that $f_\lam^{(r)}$ is
identically 0 outside $(-\lambda, \lambda)$ and then observe that
for $|x| < \lambda$,
 \[
 |f_\lam^{(r)}(x)| = |(f * \psi_\lam)^{(r)}(x)| = |(f
 * \psi_\lam^{(r)})(x)| \leq \|\psi_\lam^{(r)}\|_\infty \cdot \int_{x-\lambda}^{x+\lambda}|f|
\leq 2 B_r \lambda^{1-r}.
\]
This completes the proof.
\end{proof}\\

We now prove~(\ref{eq:superl1}), (\ref{eq:superlevi}) and
(\ref{eq:superl2}).\\

\begin{proof}
Note that the properties of $\Delta^\lam_{\lam,t}$ imply that
\begin{equation}
\label{eqn:delta-lam-props} \Pr[\boldR \leq t - 2\lambda] \leq
\E[\Delta_{\lambda, t}^\lambda(\boldR)] \leq \Pr[\boldR \leq t +
2\lambda]
\end{equation}
holds for every random variable $\boldR$ and every $t$ and $0 <
\lambda \leq 1/2$.\\

Let us first prove~(\ref{eq:superl1}), with
\[
\eps = d\,\eta^{-rd}\;\tau^{r/2 - 1}
\]
since we assume $\hone$.  Taking $\TestFcn =
\ell_1^\lam$ in Theorem~\ref{thm:cj} we obtain
 \begin{multline*}
\Bigl|\E\bigl[\ell_1(Q(\CalX))\bigr]-E\bigl[\ell_1(Q(\CalY))\bigr]\Bigr|
\leq
\Bigl|\E\bigl[\ell_1^{\lam}(Q(\CalX))\bigr]-E\bigl[\ell_1^{\lam}(Q(\CalY))\bigl]\Bigl| + 4 \lam \\
\leq (4 B_r\,\lam^{1-r} / r!)\,d\,\eta^{-rd}\;\tau^{r/2-1} + 4
\lam = O(\eps\,\lam^{1-r}) + 4\lam.
\end{multline*}
Taking $\lam = \eps^{1/r}$, gives the bound~(\ref{eq:superl1}).
Next, using~(\ref{eqn:delta-lam-props}) and applying
Theorem~\ref{thm:cj} with $\TestFcn = \Delta^\lambda_{\lambda, t}$
we obtain
\begin{multline*}
d_L(Q(\CalX), Q(\CalY)) \leq \max\left\{4\lambda, \sup_t
\Bigl|\E\bigl[\Delta^{\lam}_{\lam,t}(Q(\CalX))\bigr] -
\E\bigl[\Delta^{\lam}_{\lam,t}(Q(\CalY))\bigr]\Bigr|\right\} \\
 \leq \max\left\{(B_r\,\lambda^{-r}/r!)\,d\,\eta^{-rd}\;\tau^{r/2-1}, 4 \lam
\right\} = \max\{O(\eps\,\lam^{-r}), 4\lam\}.
\end{multline*}
Again taking $\lam = \eps^{1/(r+1)}$ we
achieve~(\ref{eq:superlevi}).  Finally, using $\TestFcn =
\trunc^\lam$ we get
 \begin{multline*}
\Bigl| \E\bigl[\trunc(Q(\CalX))\bigr] -
\E\bigl[\trunc(Q(\CalY))\bigr] \Bigr| \leq \Bigl|
\E\bigl[\trunc^{\lam}(Q(\CalX))\bigr] - \E\bigl[\trunc^{\lam}(Q(\CalY))\bigr] \Bigr| + 4 \lam^2\\
\leq (2B_{r-1}\,\lam^{2-r} / r!)\,d\,\eta^{-rd}\;\tau^{r/2-1} + 4
\lam^2 = O(\eps\,\lam^{2-r}) + 4\lam^2,
\end{multline*}
and taking $\lam = \eps^{1/r}$ we get~(\ref{eq:superl2}). This
concludes the proof of the first three bounds in
Theorem~\ref{thm:supertheorem}.
\end{proof}

\subsubsection{Closeness in distribution} \label{cid} We proceed to
prove~(\ref{eq:lim_dist}) from Theorem \ref{thm:supertheorem}.  By
losing constant factors it will suffice to prove the bound in the
case that $\CalY = \CalG$, the sequence of independent Gaussian
ensembles.  As mentioned, we will use the fact that bounded-degree
multilinear polynomials over $\CalG$ have low ``small ball
probabilities''.  Specifically, the following theorem is an
immediate consequence of Theorem~8 in~\cite{CarberyWright:01}
(taking $q=2d$ in their notation):
\begin{theorem} \label{thm:smallball}
There exists a universal constant $C$ such that for all
multilinear polynomials $Q$ of degree $d$ over $\CalG$ and all
$\eps > 0$,
\[
\P[|Q(\CalG)| \leq \eps] \leq C\,d\,(\eps/\|Q(\CalG)\|_2)^{1/d}.
\]
\end{theorem}
Thus we have the following:
\begin{corollary} \label{cor:smallball}
For all multilinear polynomials $Q$ of degree $d$ over $\CalG$
with $\Var[Q] = 1$ and for all $t \in \R$ and $\eps > 0$,
\[
\P[|Q(\CalG) - t| \leq \eps] \leq O(d\,\eps^{1/d}).
\]
\end{corollary}

We now prove~(\ref{eq:lim_dist}).\\

\begin{proof}
We will use Theorem~\ref{thm:cj} with $\TestFcn =
\Delta_{\lam,t}^{\lam}$, where $\lambda$ will be chosen later.
Writing $\Delta_t =  \Delta_{\lam,t}^{\lam}$ for brevity and using
fact~(\ref{eqn:delta-lam-props}) twice, we have
\begin{eqnarray}
\P[Q(\CalX) \leq t] & \leq & \E[\Delta_{t+2\lambda}(\CalX)] \nonumber\\
& \leq & \E[\Delta_{t+2\lambda}(\CalG)] +
|\E[\Delta_{t+2\lambda}(\CalX)] - \E[\Delta_{t+2\lambda}(\CalG)]| \nonumber\\
& \leq & \P[Q(\CalG) \leq t + 4\lam] +
|\E[\Delta_{t+2\lambda}(\CalX)] - \E[\Delta_{t+2\lambda}(\CalG)]| \nonumber\\
& = & \P[Q(\CalG) \leq t] + \P[t < Q(\CalG) \leq t+4\lam] +
|\E[\Delta_{t+2\lambda}(\CalX)] - \E[\Delta_{t+2\lambda}(\CalG)]|.
\label{eqn:d1}
\end{eqnarray}
The second quantity in~(\ref{eqn:d1}) is at most
$O(d\,(4\lam)^{1/d})$ by Corollary~\ref{cor:smallball}; the third
quantity in~(\ref{eqn:d1}) is at most $O(\eps\,\lam^{-r})$ by
Lemma~\ref{lem:approx-functional} and Theorem~\ref{thm:cj}. Thus
we conclude
\[
\P[Q(\CalX) \leq t] \leq \P[Q(\CalG) \leq t] + O(d\,\lam^{1/d}) +
O(\eps\,\lam^{-r}),
\]
independently of $t$. Similarly it follows that
\[
\P[Q(\CalX) \leq t] \geq \P[Q(\CalG) \leq t] - O(d\,\lam^{1/d}) -
O(\eps\,\lam^{-r}).
\]
independently of $t$. Choosing $\lam = \eps^{d/(rd+1)}$ we get
\[
\Bigl|\P\bigl[Q(\CalX) \leq t\bigr] - \P\bigl[Q(\CalG) \leq
t\bigr]\Bigr| \leq O(d\,\eps^{1/(rd+1)}),
\]
as required.
\end{proof}

\bigskip

The proof of Theorem \ref{thm:supertheorem} is now complete.

\subsubsection{Invariance principle for smoothed functions}
\newcommand{\Rhigh}{H(\boldR)}
\newcommand{\Rlow}{L(\boldR)}
\newcommand{\Shigh}{H(\boldS)}
\newcommand{\Slow}{L(\boldS)}
The proof of Theorem~\ref{thm:smooththeorem} is by truncating at
degree $d = c \log(1/\tau) / \log(1/\eta)$, where $c > 0$ is a
sufficiently small constant to be chosen later. Let $\Rlow =
(T_{1-\gamma} Q)^{\leq d}(\CalX)$, $\Rhigh = (T_{1-\gamma} Q)^{>
d}(\CalX)$, and define $\Slow$, and $\Shigh$ analogously for
$\CalY$.  Note that the low-degree influences of $T_{1-\gamma} Q$
are no more than those of $Q$.\\

We first prove the upper bound on $d_L(\boldR, \boldS)$.  By
Theorem \ref{thm:supertheorem} we have
\begin{equation} \label{eqn:low}
d_L(\Rlow, \Slow) \leq
d^{\Theta(1)}\,\eta^{-\Theta(d)}\,\tau^{\Theta(1)} =
\eta^{-\Theta(d)}\,\tau^{\Theta(1)}.
\end{equation}
As for $\Rhigh$ and $\Shigh$ we have $\E[\Rhigh] = \E[\Shigh] = 0$
and $\E[\Rhigh^2] = \E[\Shigh^2] \leq (1-\gamma)^{2d}$ (since
$\Var[Q] \leq 1$).  Thus by Chebyshev's inequality it follows that
for all $\lam$,
\begin{equation} \label{eqn:high}
\P[|\Rhigh| \geq \lam] \leq (1-\gamma)^{2d}/\lam^2, \qquad
\P[|\Shigh| \geq \lam] \leq (1-\gamma)^{2d}/\lam^2.
\end{equation}
Combining~(\ref{eqn:low}) and~(\ref{eqn:high}) and taking $\lam =
(1-\gamma)^{2d/3}$ we conclude that the L\'{e}vy distance between
$\boldR$ and $\boldS$ is at most
\begin{equation} \label{eqn:mainbound}
\eta^{-\Theta(d)}\,\tau^{\Theta(1)} + 4 (1-\gamma)^{2d/3} \leq
\eta^{-\Theta(d)}\,\tau^{\Theta(1)} +
\exp\bigl(-\gamma\,\Theta(d)\bigr).
\end{equation}
Our choice of $d$, with $c$ taken sufficiently small so that the
second term above dominates, completes the proof of the upper bound on $d_L(\boldR, \boldS)$.\\

To prove the claim about $\trunc$ we need the following simple
lemma:
\begin{lemma}  For all $a, b \in \R$, $|\trunc(a+b) - \trunc(a)| \leq 2|ab| +
2b^2$.
\end{lemma}
\begin{proof}
We have
\[
|\trunc(a+b) - \trunc(a)| \leq |b| \sup_{x \in [a,a+b]}
|\trunc'(x)|.
\]
The claim follows since $\trunc'(x)  = 0$ for $|x| \leq 1$ and
$|\trunc'(x)| = 2||x|-1| \leq 2|x|$ for $|x| \geq 1$.
\end{proof}

\bigskip

By~(\ref{eq:superl2}) in Theorem~\ref{thm:supertheorem} we get the
upper bound of $\eta^{-\Theta(d)}\,\tau^{\Theta(1)}$ for
$|\E[\trunc(\Rlow) - \trunc(\Slow)]|$.  The Lemma above and
Cauchy-Schwartz imply
\begin{multline*}
\E\Bigl[\bigl|\trunc(\boldR)) - \trunc(\Rlow)\bigr|\Bigr] =
\E\Bigl[\bigl|\trunc(\Rlow + \Rhigh) - \trunc(\Rlow)\bigr|\Bigr]
\leq 2 \E\bigl[|\Rlow \Rhigh|\bigr] + \E\bigl[\Rhigh^2\bigr]
\\ \leq 2\sqrt{\E[\Rhigh^2]} + \E[\Rhigh^2] \leq 2(1-\gamma)^d +
(1-\gamma)^{2d} \leq \exp\bigl(-\gamma\,\Theta(d)\bigr),
\end{multline*}
and similarly for $\boldS$. Thus
\[
|\E[\trunc(\boldR)] - \E[\trunc(\boldS)]| \leq
\eta^{-\Theta(d)}\,\tau^{\Theta(1)} +
\exp\bigl(-\gamma\,\Theta(d)\bigr)
\]
as in~(\ref{eqn:mainbound}) and we get the same upper bound.

Finally, it is easy to see that the second statement of the theorem also holds 
as the only property of $\boldR$ we have used is that $\Var[Q^{> d}] \leq (1-\gamma)^{2d}$ for all $d$.

\subsection{Invariance principle under Lyapunov conditions}
Here we sketch a proof of Theorem~\ref{thm:Lyap}.

\begin{proof} (sketch) Let $\Delta:\R \to [0,1]$ be a nondecreasing smooth function with $\Delta(0)=0,$ $\Delta(1)=1$ and $A:=\sup_{x \in \R} |\Delta'''(x)| < \infty.$ Then $\sup_{x \in \R} |\Delta''(x)| \leq A/2$ and therefore for $x,y \in \R$ we have
\[
|\Delta''(x)-\Delta''(y)| \leq A^{3-q}|\Delta''(x)-\Delta''(y)|^{q-2} \leq
A^{3-q}(A|x-y|)^{q-2}=A|x-y|^{q-2}.
\]
For $s>0$ let $\Delta_{s}(x)=\Delta(x/s),$ so that
$|\Delta_{s}''(x)-\Delta_{s}''(y)| \leq A s^{-q}|x-y|^{q-2}$
for all $x,y \in \R.$ 
Let $\boldY$ and $\boldZ$ be random variables with
$\E[\boldY]=\E[\boldZ],$ $\E[\boldY^{2}]=\E[\boldZ^{2}]$
and $\E[|\boldY|^{q}], \E[|\boldZ|^{q}]<\infty.$
Then $|\E[\Delta_{s}(x+\boldY)]-\E[\Delta_{s}(x+\boldZ)]| \leq
As^{-q}(\E[|\boldY|^{q}]+\E[|\boldZ|^{q}])$ for all $x \in \R.$
Indeed, for $u \in [0,1]$ let
$\phi(u)=\E[\Delta_{s}(x+u\boldY)]-\E[\Delta_{s}(x+u\boldZ)].$
Then $\phi(0)=\phi'(0)=0$ and
\[
|\phi''(u)|=
|\E[\boldY^{2}(\Delta_{s}''(x+u\boldY)-\Delta_{s}''(x))]-
\E[\boldZ^{2}(\Delta_{s}''(x+u\boldZ)-\Delta_{s}''(x))]| \leq
As^{-q}u^{q-2}(\E[|\boldY|^{q}]+\E[|\boldZ|^{q}]),
\]
so that $|\phi(1)| \leq As^{-q}(\E[|\boldY|^{q}]+\E[|\boldZ|^{q}]).$
Now, using the above estimate and the fact that both $\CalX$
and $\CalG$ are $(2,q,\eta)$-hypercontractive with
$\eta=\frac{\beta^{-1/q}}{2\sqrt{q-1}}$ one arrives at
\[
|\E[\Delta_{s}(Q(\boldX_{1}, \ldots, \boldX_{n}))]-
\E[\Delta_{s}(Q(\boldG_{1}, \ldots, \boldG_{n}))]| \leq
O(s^{-q}\eta^{-qd}\sum_{i}(\sum_{S \ni i} c_{S}^{2})^{q/2}).
\]
Replacing $Q$ by $Q+t$ and using the arguments of subsection~\ref{cid} 
yields
\[
 \sup_{t} \bigl|\P[Q(\boldX_1, \dots, \boldX_n) \leq t] -
          \P[Q(\boldG_1, \dots, \boldG_n) \leq t] \bigr|
  \leq
  \]
  \[
  O(ds^{1/d})+
  O(s^{-q}\eta^{-qd}\sum_{i}(\sum_{S \ni i} c_{S}^{2})^{q/2}).
\]
Optimizing over $s$ ends the proof. We skip some elementary 
calculations.
\end{proof}

\section{Proofs of the conjectures} \label{sec:conj}
Our applications of the invariance principle have the following
character: We wish to study certain noise stability properties of
low-influence functions on finite product probability spaces. By
using the invariance principle for slightly smoothed functions,
Theorem~\ref{thm:smooththeorem}, we can essentially analyze the
properties in the product space of our choosing.  And as it
happens, the necessary result for Majority Is Stablest is already
known in Gaussian space~\cite{Borell:85} and the necessary result
for It Ain't Over Till It's Over is already known on the
uniform-measure discrete cube~\cite{MORSS:04}.\\

In the case of the Majority Is Stablest problem, one needs to find
a set of prescribed Gaussian measure which maximizes the probability
that the Ornstein-Uhlenbeck process (started at the Gaussian measure)
will belong to the set at times $0$ and time $t$ for some fixed time $t$.
This problem was solved by Borell
in~\cite{Borell:85} using symmetrization arguments.  It should
also be noted that the analogous result for the sphere has been
proven in more than one place, including a paper of Feige and
Schechtman~\cite{FeigeSchechtman:02}. It fact, one can deduce Borell's result 
and Majority is Stablest from the spherical result using the proximity of 
spherical and Gaussian measures in high dimensions and the invariance
principle proven here.\\   

In the case of the It Ain't Over Till It's Over problem, the
necessary result on the discrete cube $\bits^n$ was essentially
proven in the recent paper~\cite{MORSS:04} using the reverse
Bonami-Beckner inequality (which is also due to
Borell~\cite{Borell:82}). This paper did not solve the conjecture
though (nor did that paper note the relevance), even when the
conjecture is set on $\bits^n$; the reason is that reduction of
the problem to a question about $T_\rho$ already involves
transferring to a different product domain (e.g., $\{-1, 0, 1\}^n$
with biased measure) and so the invariance principle is
required.\\

Note that in both cases the necessary auxiliary result is valid
without any assumptions about low influences. This should not be
surprising in the Gaussian case, since given a multilinear
polynomial $Q$ over Gaussians it is easy to define another
multilinear polynomial $\tilde{Q}$ over Gaussians with exactly the
same distribution and arbitrarily low influences, by letting
\[
\tilde{Q}(x_{1,1},\dots,x_{1,N},\;\ldots\;,x_{n,1},\ldots,x_{n,N})
= Q\Bigl(\frac{x_{1,1} + \cdots +
x_{1,N}}{N^{1/2}},\;\ldots\;,\frac{x_{n,1} + \cdots +
x_{n,N}}{N^{1/2}}\Bigr).
\]
The fact that low influences are not required for
the the results of~\cite{MORSS:04} is perhaps more surprising.

\subsection{Noise stability in Gaussian space}
We begin by recalling some definitions and results relevant for
``Gaussian noise stability''.  Throughout this section we consider
$\R^n$ to have the standard $n$-dimensional Gaussian distribution,
and our probabilities and expectations are over this
distribution.\\

Let $U_\rho$ denote the Ornstein-Uhlenbeck operator acting on
$L^2(\R^n)$ by
\[
(U_\rho f)(x) = \Ex_{y}[f(\rho x + \sqrt{1-\rho^2}\, y)],
\]
where $y$ is a random standard $n$-dimensional Gaussian.  It is
easy to see that if $f(x)$ is expressible as a \emph{multilinear}
polynomial in its $n$ independent Gaussian inputs,
\[
f(x_1, \dots, x_n) = \sum_{S \subseteq [n]} c_S \prod_{i \in S}
x_i,
\]
then $U_\rho f$ is the following multilinear polynomial:
\[
(U_\rho f)(x_1, \dots, x_n) = \sum_{S \subseteq [n]} \rho^{|S|}
c_S \prod_{i \in S} x_i.
\]
Thus $U_\rho$ acts identically to $T_\rho$ for multilinear
polynomials $Q$ over $\CalG$, the Gaussian sequence of
ensembles.\\

Next, given any function $f : \R^n \to \R$, recall that its
\emph{(Gaussian) nonincreasing spherical rearrangement} is defined
to be the upper semicontinuous nondecreasing function $f^* : \R
\to \R$ which is equimeasurable with $f$; i.e., for all $t \in
\R$, $f^*$ satisfies $\Pr[f > t] = \Pr[f^* > t]$ under Gaussian measure.\\

We now state a result of Borell concerning the
Ornstein-Uhlenbeck operator $U_\rho$ (see also Ledoux's
Saint-Flour lecture notes~\cite{DoGrLe:96}).  Borell uses Ehrhard
symmetrization to show the following:
\begin{theorem} (Borell~\cite{Borell:85}) \label{thm:borell}
Let $f, g \in L^2(\R^n)$.
%\[
%\sup_{x \neq y} \frac{|f(x)-f(y)|}{\|x-y\|} < \infty
%\]
%and similarly for $g$.
Then for all $0 \leq \rho \leq 1$ and all
$q \geq 1$,
\[
\E[(U_\rho f)^q \cdot g] \leq \E[(U_\rho f^*)^q \cdot g^*].
\]
\end{theorem}

Borell's result is more general and is stated for Lipschitz
functions, but standard density arguments immediately imply the
validity of the statement above. One immediate consequence of the
theorem is that $\Stab_\rho(f) \leq \Stab_\rho(f^*)$, where we
define
\begin{equation} \label{eqn:gauss-stab}
\Stab_\rho(f) = \E[f \cdot U_\rho f] = \E[(U_{\sqrt{\rho}} f)^2].
\end{equation}
One can think of this quantity as the ``(Gaussian) noise stability
of $f$ at $\rho$''; again, it is compatible with our earlier
definition of $\Stab_\rho$ if $f$ is a multilinear polynomial over
$\CalG$.\\

Note that the latter equality in~(\ref{eqn:gauss-stab}) and the
fact that $U_{\sqrt{\rho}}$ is positivity-preserving and linear
imply that $\sqrt{\Stab_\rho}$ defines an $L^2$ norm on
$L^2(\R^n),$ dominated by the usual $L^{2}$ norm, so that it is a
continuous convex functional on $L^{2}(\R^{n})$.  The set of all
$[0,1]$-valued functions from $L^{2}(\R^{n})$ having the same mean
as $f$ is closed and bounded in the standard $L^{2}$ norm and one
can easily check that its extremal points are indicator functions;
hence by the Edgar-Choquet theorem (see \cite{Edgar:75}; clearly
$L^{2}(\R^{n})$ is separable and it has the Radon-Nikodym property
since it is a Hilbert space):
\[
\sqrt{\Stab_\rho}(f) \leq \sup_{\chi}\sqrt{\Stab_\rho}(\chi),
\]
where the supremum is taken over all functions $\chi:\R^{n} \to \{
0,1\}$ with $\E[\chi]=\E[f].$ Since by Borell's result
$\Stab_\rho(\chi) \leq \Stab_\rho(\chi^{*})$, we have
$\Stab_\rho(f) \leq \Stab_\rho(\chi_{\mu})$ where $\chi_{\mu}:\R
\to \{ 0,1\}$ is the indicator function of a halfline with measure
$\mu=\E[f]$.\\
%\knote{Again, the old version is hidden behind \%. I
%made several serious changes to make things clear and correct. In
%particular after thinking a while I could not see how to apply
%Krein-Milman theorem. We cannot hope for compactness - even in the
%set of extreme points one can find an infinite separated subset;
%think of the infinite sequence of non-correlated indicator
%functions, for example. Also a standard catch which is to go to
%Lipschitz functions and use some Arzela-Ascoli stuff is hopeless
%here since our extreme points are not Lipschitz. Anyway, if you do
%not find the above explanation convincing I still have a quite
%elementary complete proof, however at least three times longer. On
%the other hand using the noncompact Choquet theorem seems too
%heavy artillery for such a simple matter... Let me know what you
%think about it. Also, please correct my English if necessary.}

Let us introduce some notation:
\begin{definition} Given $\mu \in [0,1]$, define $\chi_\mu : \R \to \{0,1\}$
to be the indicator function of the interval $(-\infty, t]$, where
$t$ is chosen so that $\E[\chi_\mu] = \mu$.  Explicitly, $t =
\Phi^{-1}(\mu)$, where $\Phi$ denotes the distribution function of
a standard Gaussian.  Furthermore, define
\[
\StabThr{\rho}{\mu} = \Stab_\rho(\chi_\mu) = \Pr[\boldX \leq t,
\boldY \leq t],
\]
where $(\boldX, \boldY)$ is a two dimensional Gaussian vector with
covariance matrix
$
\left( \begin{smallmatrix}
1 & \rho \\
\rho & 1
\end{smallmatrix} \right)
$.
\end{definition}
Summarizing the above discussion, we obtain:
\begin{corollary} \label{cor:bor}
Let $f : \R^n \to [0,1]$ be a measurable
function on Gaussian space with $\E[f] = \mu$. Then for all $0
\leq \rho \leq 1$ we have $\Stab_\rho(f) \leq
\StabThr{\rho}{\mu}$.
\end{corollary}

\noindent This is the result we will use to prove the Majority Is
Stablest conjecture.  We note that in general there is no
closed form for $\StabThr{\rho}{\mu}$; however, some asymptotics
are known: For balanced functions we have
Sheppard's formula $\StabThr{\rho}{1/2} = \frac14+
\frac{1}{2\pi}\arcsin \rho$.  Some other properties of
$\StabThr{\rho}{\mu}$ are given in Appendix~\ref{app:StabThr}.

\subsection{Majority Is Stablest} \label{subsec:maj_stablest}
In this section we prove a strengthened form of the Majority Is
Stablest conjecture.  The implications of this result were
discussed in Section~\ref{sec:misc}.
\begin{theorem} \label{thm:MIST}
Let $f : \Omega_1 \times \cdots \times \Omega_n \to [0,1]$ be a
function on a finite product probability space and assume that
for each $i$ the minimum probability of any atom in $\Omega_i$ is
at least $\alpha \leq 1/2$.  Write $K = \log(1/\alpha)$.  Further
assume that there is a $0 < \tau < 1/2$ such that $\Inf_i^{\leq
\log(1/\tau)/K}(f) \leq \tau$ for all $i$. (See
Definition~\ref{def:low-degree-influence} for the definition of
low-degree influence.)  Let $\mu = \E[f]$.  Then for any $0 \leq
\rho < 1$,
\[
\Stab_\rho(f) \leq \StabThr{\rho}{\mu} + \eps,
\]
where
\[
\eps = O\Bigl(\frac{K}{1-\rho}\Bigr) \cdot \frac{ \log \log
(1/\tau)}{\log(1/\tau)}.
\]
\end{theorem}

\rnote{In the most basic case, with rademachers and $\mu = 1/2$,
you can probably slightly improve the error term.  This may well
be important for hardness of approximation / metric space
problems, and we should check the details.  Specifically, can one
show that for any function that is not a $(1/n^{.001},
n/2)$-junta, the noise stability at $1-O(1/\log)$ is at most $1 -
\Omega(1/\sqrt{\log n})$?}

\noindent For the reader's convenience we record here two facts
from Appendix~\ref{app:StabThr}:
\begin{eqnarray*}
\StabThr{\rho}{{\textstyle \frac12}} & = & \frac14+ \frac{1}{2\pi}\arcsin \rho \\
\StabThr{\rho}{\mu} & \sim & \mu^{2/(1+\rho)}\,(4\pi\ln
(1/\mu))^{-\rho/(1+\rho)}\,\frac{(1+\rho)^{3/2}}{(1-\rho)^{1/2}}
\qquad \text{as $\mu \to 0$.}
\end{eqnarray*}
\begin{proof}
As discussed in Section~\ref{sec:setup}, let $\CalX$ be the
sequence of ensembles such that $\CalX_i$ spans the functions on
$\Omega_i$, and express $f$ as the multilinear polynomial $Q$.
We use the invariance principle under hypothesis $\htwo$.
Express
$\rho = \rho' \cdot (1-\gamma)^2$, where $0 < \gamma \ll 1-\rho$
will be chosen later. Writing $Q(x) = \sum c_{\boldsigma}
x_{\boldsigma}$ (with $c_{\mathbf{0}} = \mu$) we see that
\[
\Stab_{\rho}(Q(\CalX))  = \sum (\rho' \cdot
(1-\gamma)^2)^{|\boldsigma|} c_{\boldsigma}^2 =
\Stab_{\rho'}((T_{1-\gamma} Q)(\CalG)),
\]
where $\CalG$ is the sequence of independent Gaussian ensembles.\\

Since $Q(\CalX)$ is bounded in $[0,1]$ the same is true of $\boldR
= (T_{1-\gamma} Q)(\CalX)$. In other words, $\E[\trunc(\boldR)] =
0$, where $\trunc$ is the function from~(\ref{eq:def_trunc}).
Writing $\boldS = (T_{1-\gamma} Q)(\CalG)$, we conclude from
Theorem~\ref{thm:smooththeorem} that $\E[\trunc(\boldS)] \leq
\tau^{\Omega(\gamma/K)}$.  That is, $\|\boldS - \boldS'\|_2^2 \leq
\tau^{\Omega(\gamma/K)}$, where $\boldS'$ is the random variable
depending on $\boldS$ defined by
\[
\boldS' = \left\{ \begin{array}{rl}
                0    & \text{if $\boldS \leq 0$,} \\
                \boldS & \text{if $\boldS \in [0,1]$,} \\
                1    & \text{if $\boldS \geq 1$.}
                \end{array}
                \right.
\]
Then 
\begin{eqnarray*}
|\Stab_{\rho'}(\boldS) - \Stab_{\rho'}(\boldS')| &=& 
|\E[\boldS \cdot U_{\rho'}  \boldS] -  \E[\boldS' \cdot U_{\rho'} \boldS']| \\ 
 &\leq&   
|\E[\boldS \cdot U_{\rho'}  \boldS] -  \E[\boldS' \cdot U_{\rho'} \boldS]| + 
|\E[\boldS' \cdot U_{\rho'} \boldS] -  \E[\boldS' \cdot U_{\rho'} \boldS']| \\ 
&\leq& 
(\|\boldS\|_2 + \|\boldS'\|_2)\|\boldS-\boldS'\|_2 \leq \tau^{\Omega(\gamma/K)}, 
\end{eqnarray*} 
where we have used the fact that $U_{\rho'}$ is a contraction on $L^2$. 

Writing $\mu' = \E[\boldS']$ it follows from Cauchy-Schwartz that 
$|\mu-\mu'| \leq \tau^{\Omega(\gamma/K)}$. 
Since $\boldS'$ takes values in
$[0,1]$ it follows from Corollary~\ref{cor:bor} that
$\Stab_{\rho'}(\boldS') \leq \StabThr{\rho'}{\mu'}$.  We thus
conclude
\[
\Stab_\rho(Q(\CalX)) =  \Stab_{\rho'}(\boldS) \leq
\Stab_{\rho'}(\boldS') + \tau^{\Omega(\gamma/K)} \leq
\StabThr{\rho'}{\mu'} + \tau^{\Omega(\gamma/K)}.
\]
We can now bound the difference $|\StabThr{\rho}{\mu} -
\StabThr{\rho'}{\mu'}|$ using Lemmas~\ref{lem:I1} and
Corollary~\ref{cor:I2} in Appendix~\ref{app:StabThr}.  We get a
contribution of $2|\mu - \mu'| \leq \tau^{\Omega(\gamma/K)}$ from
the difference in the $\mu$'s and a contribution of at most
$O(\gamma/(1-\rho))$ from the difference in the $\rho$'s.  Thus we
have
\[
\Stab_\rho(Q(\CalX)) \leq \StabThr{\rho}{\mu} +
\tau^{\Omega(\gamma/K)} + O(\gamma/(1-\rho)).
\]
Taking
\[
\gamma = C \cdot K \cdot \frac{\log \log (1/\tau)}{\log(1/\tau)}
\]
for some large enough constant $C$ and this gives the claimed
bound.
\end{proof}

\subsection{It Ain't Over Till It's Over}
As mentioned, our proof of the It Ain't Over Till It's Over
conjecture will use a result due essentially to~\cite{MORSS:04}:
\begin{theorem} \label{thm:coins} Let $f : \bits^n \to [0,1]$ have $\E[f] = \mu$
(with respect to uniform measure on $\bits^n$).  Then for any $0 <
\rho < 1$ and any $0 < \eps \leq 1 - \mu$ we have
\[
\Pr[T_\rho f > 1-\delta] < \eps
\]
provided
\[
\delta < \eps^{\rho^2/(1-\rho^2) + O(\kappa)},
\]
where
\[
\kappa = \frac{\sqrt{c(\mu)}}{1-\rho} \cdot
\frac{1}{\sqrt{\log(1/\eps)}}, \qquad c(\mu) = \mu
\log(e/(1-\mu)).
\]
\end{theorem}
This theorem follows from the proof of Theorem 4.1
in~\cite{MORSS:04}; for completeness we give an explicit
derivation in Appendix~\ref{app:coins}.
\begin{remark} Since the only fact about $\bits^n$ used in the proof of
Theorem~\ref{thm:coins} is the reverse Bonami-Beckner inequality,
and since this inequality also holds in Gaussian space, we
conclude that Theorem~\ref{thm:coins} also holds for
measurable functions on Gaussian space $f : \R^n \to
[0,1]$.  In this setting the result can be proved using Borell's
Corollary~\ref{cor:bor} instead of using the reverse Bonami-Beckner 
inequality.
\end{remark}

\bigskip

The first step of the proof of It Ain't Over Till It's Over is to
extend Theorem~\ref{thm:coins} to functions on arbitrary product
probability spaces.  Note that if we only want to solve the
problem for functions on $\bits^n$ with the uniform measure, this
step is unnecessary.  The proof of the extension is very similar
to the proof of Theorem~\ref{thm:MIST}. In order to state the theorem it would 
be helpful to let $u > 0$ be a constant such that 
Theorem~\ref{thm:smooththeorem} holds with the bound $\tau^{u \gamma/K}$. 

\begin{theorem} \label{thm:general-coins}
Let $f : \Omega_1 \times \cdots \times \Omega_n \to [0,1]$ be a
function on a finite product probability space and assume that
for each $i$ the minimum probability of any atom in $\Omega_i$ is
at least $\alpha \leq 1/2$.  Let $K \geq \log(1/\alpha)$.  Further
assume that there is a $\tau > 0$ such that $\Inf_i^{\leq
\log(1/\tau)/K}(f) \leq \tau$ for all $i$ (recall 
Definition~\ref{def:low-degree-influence}). 
Let $\mu = \E[f]$.  Then for any $0 < \rho
< 1$ there exists $\eps(\mu,\rho)$ such that if $0 < \eps < \eps(\mu,\rho)$ 
we have
\[
\Pr[T_\rho f > 1 - \delta] \leq \eps
\]
provided
\[
\delta < \eps^{\rho^2/(1-\rho^2) + C \kappa}, \qquad 
\tau \leq \eps^{(100K/u(1-\rho))(1 /(1-\rho)^3 + C \kappa)}
\]
where
\[
\kappa = \frac{\sqrt{c(\mu)}}{1-\rho} \cdot
\frac{1}{\sqrt{\log(1/\eps)}}, \qquad c(\mu) = \mu
\log(e/(1-\mu)) + \eps
\]
and $C > 0$ is some constant. 
\end{theorem}
\begin{proof}
Without loss of generality we assume that 
$\delta = \eps^{\rho^2/(1-\rho^2) + C \kappa}$ 
as taking a smaller $\delta$ yields a smaller tail probability.
We can also assume $\eps(\mu,\rho)<1/10.$  
Let $\CalX$ and $Q$ be as in the proof of Theorem~\ref{thm:MIST}
and this time decompose $\rho = \rho' \cdot (1-\gamma)$ where we
take $\gamma = \kappa \cdot (1-\rho)^2$.  Note that taking 
$\eps(\mu,\rho)$ sufficiently small we have $\kappa <1, \gamma < 0.1$ and
$(1-\rho)/(1-\rho') \leq 2$. 
%Note that with
%this choice of $\gamma$ and taking $C_2$ sufficiently large 
%the quantity $\tau^{\gamma/K}$ can be bounded by 
%$\eps^{O(C_2/(1-\rho)^3)} \leq \eps^{10}$.  
Let $\boldR =
(T_{1-\gamma} Q)(\CalX)$ as before, and let $\boldS =
(T_{1-\gamma} Q)(\CalY)$, where $\CalY$ denotes the Rademacher
sequence of ensembles ($\boldY_{i,0} = 1$, $\boldY_{i,1} = \pm 1$
independently and uniformly random).  Since $\E[\trunc(\boldR)] =
0$ as before, we conclude from Theorem~\ref{thm:smooththeorem}
that we have $\E[\trunc(\boldS)] \leq 
\tau^{u \gamma/K} \leq \eps^{10/(1-\rho) + 2C \kappa}$; i.e., 
\begin{equation} \label{eq:most_ugly}
\|\boldS - \boldS'\|_2^2 \leq \eps^{10/(1-\rho) + 2C \kappa}
\end{equation}
where $\boldS'$ is the truncated version of $\boldS$ as in the
proof of Theorem~\ref{thm:MIST}. Now $\boldS'$ is a function
$\bits^n \to [0,1]$ with mean $\mu'$ differing from $\mu$ by at
most $\eps^{5}$ (using Cauchy-Schwartz, as before).  
This implies that $c(\mu') \leq O(c(\mu))$. 

Furthermore, our assumed upper
bound on $\delta$ also holds with $\rho'$ in place of $\rho$. 
This is because 
\[
\frac{{\rho'}^2}{1-{\rho'}^2} - \frac{{\rho}^2}{1-{\rho}^2} = 
\frac{1}{1-\rho'^2}-\frac{1}{1-\rho^2} \leq 
(\rho'^2-\rho^2)\frac{1}{(1-\rho'^2)^2} \leq \frac{2 \gamma}{(1-\rho')^2} 
\leq \frac{8 \gamma}{(1-\rho)^2} = 8 \kappa.
\]
Thus Theorem~\ref{thm:coins} implies that if $C$ is sufficiently large then 
\[
\Pr[T_{\rho'} \boldS' > 1 - 4 \delta] < \eps/2.
\]
This, in turn implies that  
\[
\Pr[T_{\rho'} \boldS > 1 - 2 \delta] < 3\eps/4.
\]
This follows by (\ref{eq:most_ugly}) since,
\[
\Pr[T_{\rho'} \boldS > 1 - 4 \delta] - \Pr[T_{\rho'} \boldS' > 1 - 2 \delta] 
\leq \delta^{-2} \|T_{\rho'} \boldS - T_{\rho'} \boldS'\|_2^2 \leq  
     \delta^{-2} \|\boldS - \boldS'\|_2^2.
\]  
We now use Theorem~\ref{thm:smooththeorem} again, bounding the
L\'{e}vy distance of $(T_\rho Q)(\CalY)$ and $(T_\rho Q)(\CalX)$
by $\tau^{u(1-\rho)/K}$, which is 
smaller than $\delta$ and $\eps/8$.
Thus  
\[
\Pr[(T_\rho Q)(\CalX) > 1 - \delta] \leq 
\Pr[T_\rho f > 1- 2 \delta] + \eps/8 < \eps, 
\]
as needed. 
\end{proof}

\bigskip

The second step of the proof of It Ain't Over Till It's Over is to
use the invariance principle to show that the random variable
$V_\rho f$ (recall Definition~\ref{def:V}) has essentially the
same distribution as $T_{\sqrt{\rho}} f$.

\begin{theorem} \label{thm:vote-dist}
Let $0 < \rho < 1$ and let $f : \Omega_1 \times \cdots \times
\Omega_n \to [0,1]$ be a function on a finite product
probability space; assume that for each $i$ the minimum
probability of any atom in $\Omega_i$ is at least $\alpha \leq
1/2$.  Further assume that there is a $0 < \tau < 1/2$ such that
$\Inf_i^{\leq\,\log(1/\tau)/K'} \leq \tau$ for all $i$, where $K'
= \log(1/(\alpha \rho (1-\rho)))$.  Then
\[
d_L(V_\rho f, T_{\sqrt{\rho}} f) \leq \tau^{\Omega((1-\rho)/K')}.
\]
\end{theorem}
\begin{proof}
Introduce $\CalX$ and $Q$ as in the proof of
Theorems~\ref{thm:MIST} and~\ref{thm:general-coins}.  We now
define a new independent sequence of orthonormal ensembles
$\CalX^{(\rho)}$ as follows.  Let $\boldV_1, \dots, \boldV_n$ be
independent random variables, each of which is $1$ with
probability $\rho$ and $0$ with probability $1-\rho$.  Now define
$\CalX^{(\rho)} = (\CalX_1^{(\rho)}, \dots, \CalX_n^{(\rho)})$ by
$\boldX^{(\rho)}_{i, 0} = 1$ for each $i$, and
$\boldX^{(\rho)}_{i,j} = \rho^{-1/2} \boldV_i \boldX_{i,j}$ for
each $i$ and $j > 0$.  It is easy to verify that $\CalX^{(\rho)}$
is indeed an independent sequence of orthonormal ensembles.  We
will also use the fact that each atom in the ensemble
$\CalX_i^{(\rho)}$ has weight at least $\alpha' = \alpha \cdot
\min\{\rho, 1-\rho\} \geq \alpha \rho(1-\rho)$.
(one can also use Proposition~\ref{prop:add} to get a bit better
estimate on $K'$).\\

The crucial observation is now simply that the random variable
$V_\rho f$ has precisely the same distribution as the random
variable $(T_{\sqrt{\rho}} Q)(\CalX^{(\rho)})$. The reason is that
when the randomness in the $\boldV_i = 1$ ensembles is fixed, the
expectation of the restricted function is given by substituting
$0$ for all other random variables $\boldX_{i,j}$. The
$T_{\sqrt{\rho}}$ serves to cancel the factors $\rho^{-1/2}$
introduced in the definition of $\boldX^{(\rho)}_{i,j}$ to ensure orthonormality.\\

It now simply remains to use Theorem~\ref{thm:smooththeorem} to
bound the L\'{e}vy distance of  $(T_{\sqrt{\rho}}
Q)(\CalX^{(\rho)})$ and $(T_{\sqrt{\rho}} Q)(\CalX)$, where here
$\CalX$ denotes a copy of this sequence of ensembles.
We use hypothesis $\hthree$ and
get a bound of $\tau^{\Omega((1-\sqrt{\rho})/K')} =
\tau^{\Omega((1-\rho)/K')}$, as required.
\end{proof}

\bigskip

Our generalization of It Ain't Over Till It's Over is now simply a
corollary of Theorems~\ref{thm:general-coins}
and~\ref{thm:vote-dist}; by taking $K'$ instead of $K$ in the
upper bound on $\tau$ and taking $\delta$ to have its maximum possible value, 
 we make the error of
\[
\tau^{u((1-\rho)/K')} \leq 
\eps^{(100/(1-\rho))(1 /(1-\rho)^3 + C \kappa)}
\] 
from Theorem~\ref{thm:vote-dist} which is 
negligible compared to both $\eps$ and $\delta$ below.
\begin{theorem} \label{thm:aint}
Let $0 < \rho < 1$ and let $f : \Omega_1 \times \cdots \times
\Omega_n \to [0,1]$ be a function on a finite product
probability space; assume that for each $i$ the minimum
probability of any atom in $\Omega_i$ is at least $\alpha \leq
1/2$.  Further assume that there is a $0 < \tau < 1/2$ such that
$\Inf_i^{\leq\,\log(1/\tau)/K'} \leq \tau$ for all $i$, where $K'
= \log(1/(\alpha \rho (1-\rho)))$.  Let $\mu = \E[f]$. 
Then there exists an $\eps(\rho,\mu) > 0$ such that if 
$\eps \leq \eps(\rho,\mu)$ then
\[
\Pr[V_\rho f > 1 - \delta] \leq \eps
\]
provided
\[
\delta < \eps^{\rho^2/(1-\rho^2) + C \kappa}, \qquad 
\tau \leq \eps^{(100K'/u(1-\rho))(1 /(1-\rho)^3 + C \kappa)}
\]
where
\[
\kappa = \frac{\sqrt{c(\mu)}}{1-\rho} \cdot
\frac{1}{\sqrt{\log(1/\eps)}}, \qquad c(\mu) = \mu
\log(e/(1-\mu)) + \eps,
\]
where $C > 0$ is some finite constant. 
\end{theorem}
\begin{remark}
To get $V_\rho f$ bounded away from both $0$ and $1$ as desired in
Conjecture~\ref{conj:ain't}, simply use Theorem~\ref{thm:aint}
twice, once with $f$, once with $1-f$.
\end{remark}

\section{Weight at low levels --- a counterexample}
\label{sec:counterexample}

The simplest version of the Majority Is
Stablest result states roughly that among all balanced functions
$f : \{-1,1\}^n \to \{-1,1\}$ with small influences, the Majority
function maximizes $\sum_S \rho^{S} \hat{f}(S)^2$ for each $\rho$.
One might conjecture that more is true; specifically, that
Majority maximizes $\sum_{|S| \leq d} \hat{f}(S)^2$ for each $d =
1, 2, 3, \dots$.  This is known to be the case for $d = 1$
(\cite{KKMO:04}) and is somewhat suggested by the theorem of
Bourgain~\cite{Bourgain:02} which says that $\sum_{|S| \leq d}
\hat{f}(S)^2 \leq 1 - d^{-1/2 - o(1)}$ for functions with low
influences.  An essentially weaker conjecture was made
Kalai~\cite{Kalai:02}:
\begin{conjecture} \label{conj:kalai}
Let $d \geq 1$ and let $C_n$ denote the collection of all
functions $f : \{-1,1\}^n \to \{-1,1\}$ which are odd and
transitive-symmetric (see Section~\ref{sec:misc-discuss}'s
discussion of~\cite{Kalai:02}). Then
\[
\limsup_{n \to \infty} \sup_{f \in C_n} \sum_{|S| \leq d}
\hat{f}(S)^2 = \lim_{\text{$n$ odd } \to \infty} \sum_{|S| \leq d}
\widehat{\Maj_n}(S)^2.
\]
\end{conjecture}

We now show that these conjectures are false:  We construct a
sequence $(f_n)$ of completely symmetric odd functions with small
influences that have more weight on levels $1$, $2$, and $3$ than
Majority has.  \ignore{By odd we mean that $f(-x) = -f(x)$;} By
``completely symmetric'' we mean that $f_n(x)$ depends only on
$\sum_{i=1}^n x_i$; because of this symmetry our counterexample is
more naturally viewed in terms of the Hermite expansions of
functions $f : \R \to \{-1,1\}$ on one-dimensional Gaussian
space.\\

There are several normalizations of the Hermite polynomials in the
literature.  We will follow~\cite{LedouxTalagrand:91} and define
them to be the orthonormal polynomials with respect to the
one-dimensional Gaussian density function $\phi(x) =
e^{-x^2/2}/\sqrt{2\pi}$. Specifically, we define the Hermite
polynomials $h_d(x)$ for $d \in \N$ by
\[
\exp(\lambda x - \lambda^2/2) = \sum_{d=0}^\infty
\frac{\lambda^d}{\sqrt{d!}}\;h_d(x).
\]
The first few such polynomials are $h_0(x) = 1$, $h_1(x) = x$,
$h_2(x) = (x^2 - 1)/\sqrt{2}$, and $h_3(x) = (x^3 - 3x)/\sqrt{6}$.
The orthonormality condition these polynomials satisfy is
\[
\int_\R h_d(x) h_{d'}(x) \phi(x)\,dx = \left\{\begin{array}{cl} 1
& \text{if $d = d'$,} \\ 0 & \text{else.} \end{array} \right.
\]

\ignore{
 More concretely, we will use the following easy fact.
\begin{lemma} \label{lem:hermite}
Let $f : \R \to \R$ be a bounded Riemann measurable function.
Define $f_n : \{-1,1\}^n \to \R$ by letting
\[
f_n(x_1,\ldots,x_n) = f \left(\frac{1}{\sqrt{n}} \sum_{i=1}^n x_i \right).
\]
Then for all $d$ it holds that
\[
\lim_{n \to \infty} \sum_{|S|=d} \hat{f_n}^2(S) =
\int_{\R} h_d(x) f(x) \phi(x) dx.
\]
\end{lemma}
}

We will actually henceforth consider functions whose domain is
$\R^* = \R \setminus \{0\}$, for simplicity; the value of a
function at a single point makes no difference to its Hermite
expansion.  Given a function $f : \R^* \to \R$ we write
$\hat{f}(d)$ for $\int h_d(x) f(x) \phi(x)\,dx$. Let us also use
the notation $\Maj$ for the function which is $1$ on $(0, \infty)$
and $-1$ on $(-\infty, 0)$.

\begin{theorem}  \label{thm:counter} There is an odd function $f : \R^* \to \{-1,1\}$ with
\[
\sum_{d \leq 3} \hat{f}(d)^2 \geq .75913 > \frac{2}{\pi} +
\frac{1}{3\pi} = \sum_{d \leq 3} \widehat{\Maj}(d)^2.
\]
\end{theorem}
\begin{proof}
Let $t > 0$ be a parameter to be chosen later, and let $f$ be the
function which is $1$ on $(-\infty, -t]$ and $(0, t)$, and $-1$ on
$(-t, 0)$ and $[t, \infty)$.  Since $f$ is odd, $\hat{f}(0) =
\hat{f}(2) = 0$.   Elementary integration gives
\[
F_1(t) = \int h_1(x) \phi(x)\,dx = -e^{-t^2/2}/\sqrt{2\pi}, \qquad
F_3(t) = \int h_3(x) \phi(x)\,dx = (1-t^2)e^{-t^2/2}/\sqrt{12\pi};
\]
thus
\begin{eqnarray*}
\hat{f}(1) & = & 2(F_1(t) + F_1(-t) - F_1(0)) - F_1(\infty) -
F_1(-\infty) = \sqrt{2/\pi}\,(1-2e^{-t^2/2}), \\
\hat{f}(3) & = & 2(F_1(t) + F_1(-t) - F_1(0)) - F_1(\infty) -
F_1(-\infty) = -\sqrt{1/3\pi}\,(1 - 2(1-t^2)e^{-t^2/2}).
\end{eqnarray*}
We conclude
\begin{equation} \label{eqn:formula}
\sum_{d \leq 3} \hat{f}(d)^2 = \frac{2}{\pi} \Bigl(1 -
2e^{-t^2/2}\Bigr)^2 + \frac{1}{3\pi} \Bigl(1 -
2(1-t^2)e^{-t^2/2}\Bigr)^2.
\end{equation}
As $t \to 0$ or $\infty$ we recover the fact, well known in the
boolean regime (see, e.g., \cite{Bernasconi:98}), that $\sum_{d
\leq 3} \widehat{\Maj}(d)^2 = 2/\pi + 1/3\pi$. But the above
expression is not maximized for these $t$; rather, it is maximized
at $t = 2.69647$, where the expression becomes roughly $.75913$.
Fixing this particular $t$ completes the proof.
\end{proof}

\bigskip

It is now clear how to construct the sequence of completely
symmetric odd functions $f_n : \bits^n \to \bits$ with the same
property --- take $f_n(x) = f((x_1 + \cdots + x_n)/\sqrt{n})$. The
proof that the property holds follows essentially from the fact
that the limits of Kravchuk polynomials are Hermite polynomials.
For completeness, give a direct proof of
Corollary~\ref{cor:counter} in Appendix~\ref{app:counter}.

\begin{corollary} \label{cor:counter}
For $n$ odd there is a sequence of completely symmetric odd
functions $f_n : \{-1,1\}^n \to \{-1,1\}$ satisfying $\Inf_i(f_n)
\leq O(1/\sqrt{n})$ for each $i$, and
\[
\lim_{n \mbox{ odd }\to \infty} \sum_{|S| \leq 3}
\widehat{f_n}(S)^2 \geq 0.75913 > \frac{2}{\pi} + \frac{1}{3\pi} =
\lim_{n \mbox{ odd }\to \infty} \sum_{|S| \leq 3}
\widehat{\Maj_n}(S)^2.
\]
\end{corollary}

\bigskip

In light of this counterexample, it seems we can only hope to
sharpen Bourgain's Theorem~\ref{thm:bourgain} in the asymptotic
setting; one might ask whether its upper bound can be improved to
\[
1 - (1-o(1)) (2/\pi)^{3/2}\;d^{-1/2},
\]
the asymptotics for Majority.

%\section{Conclusion}
%In this paper we have given a new framework for analyzing boolean
%functions, and more generally, functions on discrete product
%probability spaces.  We have identified influences as a crucial
%parameter of such functions and shown that bounded-degree
%functions with low influences act like functions in Gaussian
%space.  We have further shown that \emph{all} functions with low
%influences act like functions in Gaussian space if they are
%slightly smoothed by an application of $T_{1-\gamma}$; this
%allowed us to prove several new results about the noise stability
%of functions on discrete product spaces.  We strongly believe our
%new techniques will find further application in the harmonic
%analysis of boolean functions, in social choice, in PCP analysis,
%and in combinatorics.
%\rnote{Anything else we want to say?}

\bibliographystyle{abbrv}
\bibliography{all,my}

\appendix

\section{Hypercontractivity of sequences of ensembles} \label{app:hc}

We now give the proofs of Propositions~\ref{prop:join-hypercon}
and~\ref{prop:hypercon} that were omitted.  As mentioned, the
proofs are completely straightforward adaptations of the analogous
proofs in~\cite{Janson:97}. \\

\begin{proof} (of Proposition~\ref{prop:join-hypercon})\ \
Let $Q$ be a multilinear polynomial over $\CalX \cup \CalY$.  Note
that we can write $Q(\CalX \cup \CalY)$ as $ \sum_j c_j
\CalX_{\boldsigma_j} \CalY_{\boldupsilon_j}, $ where the
$\boldsigma$'s are multi-indexes for $\CalX$, the $\boldupsilon$'s
are multi-indexes for $\CalY$, and the $c_j$'s are constants. Then
\begin{eqnarray*}
\| (T_\eta Q)(\CalX \cup \CalY) \|_q & = & \|\sum_j
\eta^{|\sigma_j| + |\upsilon_j|}  c_j \CalX_{\sigma_j}
\CalY_{\upsilon_j} \|_q \\
& = & \| \| T_\eta^{(\CalY)} \Bigl( \sum_j (\eta^{|\sigma_j|} c_j
\CalX_{\sigma_j}) \CalY_{\upsilon_j} \Bigr)
\|_{L^q(\CalY)} \|_{L^q(\CalX)} \\
& \leq & \| \| \sum_j (\eta^{|\sigma_j|} c_j \CalX_{\sigma_j})
\CalY_{\upsilon_j}
\|_{L^p(\CalY)} \|_{L^q(\CalX)} \\
& \leq & \| \| \sum_j (\eta^{|\sigma_j|} c_j \CalX_{\sigma_j})
\CalY_{\upsilon_j}
\|_{L^q(\CalX)} \|_{L^p(\CalY)} \\
& = & \| \| T_\eta^{(\CalX)} \Bigr(\sum_j (c_j \CalY_{\upsilon_j})
\CalX_{\sigma_j} \Bigl) \|_{L^q(\CalX)}
\|_{L^p(\CalY)} \\
& \leq & \| \sum_j c_j \CalY_{\upsilon_j} \CalX_{\sigma_j}
\|_{L^p(\CalX)} \|_{L^p(\CalY)} \\
& = & \|Q(\CalX \cup \CalY)\|_p,
\end{eqnarray*}
where the second inequality used a simple consequence of the
integral version of Minkowski's inequality and $p \leq q$
(see~\cite[Prop.\ C.4]{Janson:97}).
\end{proof}

\bigskip

\begin{proof} (of Proposition~\ref{prop:hypercon})\ \
Note that if $Q = Q^{=d}$ then we obviously have equality.  In the
general case, write $Q = \sum_{i = 0}^d Q^{= i}$, and note that
$\E[Q^{=i}(\CalX) Q^{=j}(\CalX)] = 0$ for $i \neq j$ is easy to
check. Thus
\begin{eqnarray*}
\|Q(\CalX)\|_q &=& \|T_{\eta} \Bigl( \sum_{i=0}^d \eta^{-i}
Q^{=i}(\CalX) \Bigr) \|_q\;\leq\;\|\sum_{i=0}^d \eta^{-i}
Q^{=i}(\CalX) \|_2 \\ &=& \left( \sum_{i=0}^d \eta^{-2i} \|
Q^{=i}(\CalX) \|_2^2 \right)^{1/2} \leq\;\eta^{-d} \| Q(\CalX)
\|_2.
\end{eqnarray*}
\end{proof}

\medskip

Let us also mention some standard facts about the
$(2,q,\eta)$-hypercontractivity of random variables. Let $q>2$.
Clearly, if we want $\boldX$ to be $(2,q,\eta)$-hypercontractive,
we must assume that $\E[|\boldX|^{q}] < \infty.$ If $\boldX$ is
$(2,q,\eta)$-hypercontractive for some $\eta \in (0,1)$ then
$\E[\boldX]=0$ and $\eta \leq (q-1)^{-1/2}.$ Indeed, it suffices
to consider the first and second order Taylor expansions in both
sides of the inequality $\| 1+\eta b\boldX \|_{q} \leq \|
1+b\boldX \|_{2}$ as $b \to 0$.  We leave details to the reader.\\

We now give the proofs of Proposition~\ref{prop:bdd} 
and Proposition~\ref{prop:add}, which follow
the argument of Szulga~\cite[Prop. 2.20]{Szulga:98}:\\

\begin{proof} (of Proposition~\ref{prop:bdd})\ \
Let $\boldX'$ be an independent copy of $\boldX$ and put
$\boldY=\boldX-\boldX'.$ 
\knote{The next sentence is added}
By the triangle inequality $\| \boldY \|_{q} \leq 2\| \boldX \|_{q}.$
Let $\epsilon$ be a symmetric $\pm 1$
Bernoulli random variable independent of $\boldY.$ Note that
$\boldY$ is symmetric, so it has the same distribution as
$\epsilon \boldY$.  Now by Jensen's inequality, Fubini's theorem,
$(2,q,(q-1)^{-1/2})$-hypercontractivity of $\epsilon$, and
Minkowski's inequality we get
\[
\| a+\eta_{q}\boldX \|_{q} \leq \| a+\eta_{q}\boldY \|_{q} = \|
a+\eta_{q}\epsilon \boldY \|_{q} \leq (\E_{\bold Y}[(\E_{\epsilon}
[|a+(q-1)^{1/2}\eta_{q} \epsilon \boldY|^{2}])^{q/2}])^{1/q}=
\]
\[
(\E[(a^{2}+(q-1)\eta_{q}^{2}\boldY^{2})^{q/2}])^{1/q}= \|
a^{2}+(q-1)\eta_{q}^{2}\boldY^{2} \|_{q/2}^{1/2} \leq
(a^{2}+(q-1)\eta_{q}^{2}\| \boldY^{2} \|_{q/2})^{1/2}=
\]
\[
(a^{2}+(\frac{\| \boldY \|_{q}}{2\| \boldX \|_{q}})^{2} \cdot \E[\boldX^{2}])^{1/2} \leq
(a^{2}+\E[\boldX^{2}])^{1/2}= \|
a+\boldX \|_{2}.
\]
\end{proof}

\bigskip

\begin{proof} (of Proposition~\ref{prop:add})\ \
Let $(\boldX',\boldV')$ be an independent copy of $(\boldX, \boldV)$
and put $\boldY=\boldV \boldX - \boldV' \boldX'.$ Then
$\| \boldY \|_{q} \leq 2\| \boldV \|_{q} \| \boldX \|_{q}=
2\rho^{1/q}\| \boldX \|_{q}$ and as in the previous proof
we get
\[
\| a+\xi_{q}\boldV \boldX \|_{q} \leq
\| a+\xi_{q}\boldY \|_{q} \leq
(a^{2}+(q-1)\xi_{q}^{2}\| \boldY \|_{q}^{2})^{1/2} \leq
(a^{2}+4(q-1)\xi_{q}^{2}\rho^{2/q}\| \boldX \|_{q}^{2})^{1/2}=
\| a+ \boldV \boldX \|_{2}. 
\]
\end{proof}

\bigskip

If $\boldX$ is defined on a finite probability space in which probability of all atoms is at least $\alpha$ then obviously
$\E[\boldX^{2}] \geq \alpha \| \boldX \|_{\infty}^{2},$ so that
$\E[|\boldX |^{q}] \leq \|
\boldX \|_{\infty}^{q-2} \cdot \E[\boldX^{2}] \leq
(\E[\boldX^{2}])^{q/2} \alpha^{1-\frac{q}{2}}.$ In particular,
if $q=3$ then $\| \boldX \|_{3}/\| \boldX \|_{2} \leq \alpha^{-1/6},$
so that $\boldV \boldX$ is $(2,3,\xi_{3})$-hypercontractive
with $\xi_{3}=2^{-3/2}\alpha^{1/6}\rho^{1/6}.$
\knote{Here I couldn't decide: shall we make changes in the paper, wherever K' appears and so on, or leave the things as they are now,
just indicating that the above observation together with the above
proof give better asymptotics for K' when rho is close to one?}

\bigskip

Let us also point out that if $\E[\boldX^{4}]<\infty$ and $\boldX$
is symmetric then a direct and elementary calculation shows that
$\boldX$ is $(2,4,\eta_{4})$-hypercontractive with
$\eta_{4}=\min(3^{-1/2}, \| \boldX \|_{2}/\| \boldX \|_{4})$ and
the constant is optimal. Therefore the random variable
$\boldX_{i,j}^{(\rho)}$ which appears in the proof of
Theorem~\ref{thm:vote-dist} is
$(2,4,\min(\rho^{1/4},3^{-1/2}))$-hypercontractive if $\CalX$ is
the $\pm 1$ Rademacher ensemble; this can be used to get a smaller
value for $K'$ if $\rho$ is close to 1. \knote{Ryan, I think this
is really much better than what is in the general case where
factor $(1-\rho)$ appears, which I guess must spoil a lot when
$\rho$ goes to $1$.}

\section{Properties of $\StabThr{\rho}{\mu}$} \label{app:StabThr}

\rnote{I meant to include some basic properties here, but I got
lazy} \ignore{ We start with some simple facts about
$\StabThr{\rho}{\mu}$.
\begin{proposition}
\begin{enumerate}
 \item $\StabThr{\rho}{1-\mu} = 1 - 2\mu + \StabThr{\rho}{\mu}$ and thus it suffices to understand $\StabThr{\rho}{\mu}$ for
$\mu \leq 1/2$.
 \item $\StabThr{\rho}{1/2 + \mu'/2} = 1/4 + \mu'/2 + \Stab_\rho(
\end{enumerate}
\end{proposition}
}

Sheppard's Formula~\cite{Sheppard:99} gives the value of
$\StabThr{\rho}{1/2}$:
\begin{theorem} \label{thm:sheppard}
$\StabThr{\rho}{1/2} = \frac14+ \frac{1}{2\pi}\arcsin \rho$.
\end{theorem}

For fixed $\rho$, the asymptotics of $\StabThr{\rho}{\mu}$ as $\mu
\to 0$ can be determined precisely; calculations of this nature
appear in~\cite{RinotRotar:01,dKPaWa:u}.
\begin{theorem} As $\mu \to 0$,
\[
 \StabThr{\rho}{\mu} \sim \mu^{2/(1+\rho)}\,(4\pi\ln(1/\mu))^{-\rho/(1+\rho)}\,\frac{(1+\rho)^{3/2}}{(1-\rho)^{1/2}}.
\]
\end{theorem}
\begin{proof}
This follows from, e.g., Lemma 11.1 of~\cite{dKPaWa:u}; although we
have $\rho > 0$ as opposed to $\rho < 0$ as in~\cite{dKPaWa:u}, the
formula there can still be seen to hold when $x = y$ (in their
notation).
\end{proof}
\rnote{How small $\mu$ has to be as a function of $1-\rho$ before
the above asymptotics kicks in I don't know.}

\begin{lemma} \label{lem:I1} For all $0 \leq \rho \leq 1$ and all $0 \leq \mu_1 \leq \mu_2 \leq 1$,
\[
\StabThr{\rho}{\mu_2} - \StabThr{\rho}{\mu_1} \leq 2(\mu_2 -
\mu_1).
\]
\end{lemma}
\begin{proof}
Let $\boldX$ and $\boldY$ be $\rho$-correlated Gaussians and write
$t_i = \Phi^{-1}(\mu_i)$. Then
\begin{eqnarray*}
\StabThr{\rho}{\mu_2} - \StabThr{\rho}{\mu_1} & = & \Pr[\boldX
\leq t_2,
\boldY \leq t_2] - \Pr[\boldX \leq t_1, \boldY \leq t_1] \\
& \leq & 2\Pr[t_1 \leq \boldX \leq t_2]\;\;=\;\;2(\mu_2
-\mu_1).
\end{eqnarray*}
\end{proof}

\begin{lemma} \label{lem:I2} Let $0 < \mu < 1$ and $0 < \rho_1 < \rho_2 <
1$, and write $I_2 = (\StabThr{\rho_2}{\mu} - \mu^2)/\rho_2$. Then 
$I_2 \leq \mu$ and 
\[
\StabThr{\rho_2}{\mu} - \StabThr{\rho_1}{\mu} \leq 4\cdot
\frac{1+\ln(\mu/I_2)}{1-\rho_2} \cdot I_2 \cdot \left(\rho_2 -
\rho_1\right).
\]
\end{lemma}
\begin{proof}
Let
\[
d = \frac{1+\ln(\mu/I_2)}{1-\rho_2}.
\]
The proof will rely on the fact that $\Gamma_{\rho}{(\mu)}$ is a convex 
function of $\rho$. This implies in particular that $I_2 \leq \mu$. 
Moreover, by the Mean Value Theorem it
suffices to show that the derivative at $\rho_2$ is at most $4 d
I_2$.  If we write the Hermite polynomial expansion of $\chi_\mu$
as $\chi_\mu(x) = \sum_{i} c_i H_i(x)$,
then $\StabThr{\rho}{\mu} = \sum_{i} c_i^2
\rho^{i}$, and thus
 \begin{equation}
 \label{eqn:d_rho} \left.
\frac{\partial}{\partial \rho} \StabThr{\rho}{\mu} \right|_{\rho =
\rho_2} = \sum_{i \geq 1} i c_i^2 \rho_2^{i - 1} \leq 
\sum_{1 \leq i \leq d+1} i c_i^2 \rho_2^{i - 1} +
\sum_{i \geq d+1} i c_i^2 \rho_2^{i - 1}.
 \end{equation}
 We will complete the proof by showing that both terms
 in~(\ref{eqn:d_rho}) are at most $2d I_2$.  
The first term is visibly at most $(d+1) I_2
\leq 2d I_2$. As for the second term, the quantity $i \rho_2^{i-1}$ 
decreases for $i \geq \rho_2/(1-\rho_2)$.  
Since $d+1 \geq (2-\rho)/(1-\rho) \geq
\rho_2/(1-\rho_2)$ the second term is therefore at most $(d+1)
\rho_2^{d} I_2 \leq (d+1) \rho_2^d \mu$.  But
\[
\rho_2^{d} \leq \rho_2^{\ln(\mu/I_2) / (1-\rho_2)} \leq I_2/\mu
\]
since $\rho_2^{1/(1-\rho_2)} \leq 1/e$ for all $\rho_2$. Thus the
second term of~(\ref{eqn:d_rho}) is also at most $(d+1) I_2 \leq
2d I_2$, as needed.
\end{proof}

Using the fact that $-I_2 \ln I_2$ is a bounded quantity we obtain: 
\begin{corollary}
\label{cor:I2} For all $0 < \mu < 1$ and $0 < \rho < 1$, if $0 <
\delta < (1-\rho)/2$ then
\[
\StabThr{\rho + \delta}{\mu} - \StabThr{\rho}{\mu} \leq
\frac{O(1)}{1-\rho} \cdot \delta.
\]
\end{corollary}

\section{Proof of Theorem~\ref{thm:coins}} \label{app:coins}

\begin{proof}  The proof is essentially the same as the proof of
the ``upper bound'' part of the proof of Theorem~4.1
in~\cite{MORSS:04}.  By way of contradiction, suppose the upper
bound on $\delta$ holds and yet $\P[T_\rho f > 1-\delta\}] \geq \eps$. 
Let $g$ be the indicator function of a subset of
$\{x : T_\rho f(x) > 1-\delta\}$ 
whose measure is $\eps$.     
%Note that by decreasing the value of $\delta$ 
%we may assume without loss of generality that $\E[g] = \eps$.  
%We will assume $\E[g] = \eps$, as
%this only hurts our argument.\\

Let $h = {\bf 1}_{\{f \leq b\}}$, where $b = 1/2 + \mu/2$. By a
Markov argument, 
\[
\mu = \E[f] \geq (1-\E[h])b \implies \E[h] \geq 1 - \E[f]/b = 
\frac{1-\mu}{1+\mu}.
\]
By another Markov argument, whenever $g(x) = 1$ we have 
\[
T_{\rho} (1-f) < \delta \implies 
T_{\rho}(h(1-b)) < \delta 
\implies T_{\rho}h < \frac{\delta}{1-b}.
\] 
Thus
\begin{equation} \label{eqn:contra1}
\E[g T_\rho h] \leq \frac{2 \eps \delta}{1-\mu}.
\end{equation}
But by Corollary~3.5 in~\cite{MORSS:04} (itself a simple corollary
of the reverse Bonami-Beckner inequality),
\begin{equation} \label{eqn:contra2}
\E[g T_\rho h] \geq \eps \cdot \eps^{(\sqrt{\alpha} +
\rho)^2/(1-\rho^2)},
\end{equation}
where $\alpha = \log(1/\E[h]) / \log(1/\eps)$.  (In Gaussian
space, this fact can also be proven using Borell's
Corollary~\ref{cor:bor}.)  Note that since $\E[h] \geq
(1-\mu)/(1+\mu)$ we get $\alpha \leq O(c(\mu)/\log(1/\eps))$,
which is also at most $1$ since we assume $\eps \leq 1-\mu$. Therefore the exponent 
$(\sqrt{\alpha} + \rho)^2$ is $\rho^2 + O(\sqrt{\alpha})$. 
Now~(\ref{eqn:contra2}) implies
\begin{equation} \label{eqn:contra3}
\E[g T_\rho h] \geq \eps \cdot \eps^{\rho^2/(1-\rho^2)} \cdot
\eps^{O(\sqrt{c(\mu)/\log(1/\eps)}/(1-\rho))} = 
\eps \cdot
\eps^{\rho^2/(1-\rho^2) + O(\kappa)}.
\end{equation}
Combining~(\ref{eqn:contra1}) and~(\ref{eqn:contra3}) yields
\[
\delta \geq ({\textstyle \frac{1-\mu}{2}}) \cdot
\eps^{\rho^2/(1-\rho^2) + O(\kappa)} 
= \eps^{\log(2/(1-\mu)) / \log(1/\eps)} 
\eps^{\rho^2/(1-\rho^2) + O(\kappa)} = \eps^{\rho^2/(1-\rho^2) + O(\kappa)},
\]
a contradiction.
\end{proof}

\section{Proof of Corollary~\ref{cor:counter}} \label{app:counter}
\begin{proof}
We define $f_n$ by setting $1 \leq u < n$ to be the odd integer
nearest to $t \sqrt{n}$ (where $t$ is the number chosen in
Theorem~\ref{cor:counter}) and then taking
\[
f_n(x) = \left\{ \begin{array}{rlcl}
                    1 & \text{if $|x| \in [1,u]$} &\text{or}& |x| \in [-n,-(u+2)],\\
                    -1 & \text{if $|x| \in [u+2,n]$} &\text{or}& |x| \in [-u,-1],\\
                 \end{array}
         \right.
\]
where $|x|$ denotes $\sum_{i=1}^n x_i$.  This is clearly a
completely symmetric odd function.  It is well known that for any
boolean function, $\sum_{i=1}^n \Inf_i(f_n)$ equals the expected
number of pivotal bits for $f_n$ in a random input. One can easily
see that this is $O(\sqrt{n})$. Thus each of $f_n$'s coordinates
has influence $O(1/\sqrt{n})$, by symmetry.\\

Let $p(n,s)$ denote the probability that the sum of $n$
independent $\pm 1$ Bernoulli random variables is exactly $s$, so
\[
p(n,s) = 2^{-n} {n \choose \frac12 n + \frac12 s},
\]
and for a set $S$ of integers let $p(n,S)$ denote $\sum_{s \in S} p(n,s)$.\\

By symmetry all of $f_n$'s Fourier coefficients at level $d$ have
the same value; we will write $\widehat{f_n}(d)$ for this
quantity. Since $f_n$ is odd, $\widehat{f_n}(0) = \widehat{f_n}(2)
= 0$.  By explicit calculation, we have
\[
\widehat{f_n}(1) = p(n-1,0) - 2p(n-1, u+1) %=\sqrt{\frac{2}{\pi n}} (1 - 2 e^{-u^2/2}) (1 + o(1))
\]
and
\begin{eqnarray*}
\widehat{f_n}(3) & = & \frac14 \Bigl(p(n-3, \{-2,2\}) -
2p(n-3,0)\Bigr) \\ & & \qquad \qquad - \frac14 \Bigl( p(n-3,
\{\pm(u-1), \pm(u+3)\}) - 2p(n-3, \{\pm (u+1)\})\Bigr) \\
& = & -\frac{1}{n-1}\;p(n-3,0) + 2\;\frac{(n-1) - (u+1)^2}{(n-1)^2
- (u+1)^2}\;p(n-3,u+1),
\end{eqnarray*}
where the last equality is by explicit conversion to factorials
and simplification.  Using $p(n, t\sqrt{n}) =
(1+o(1))\sqrt{2/\pi}e^{-t^2/2}\;n^{-1/2}$ as $n \to \infty$, we
conclude
\[
\widehat{f_n}(1) \sim \sqrt{2/\pi}(1-2e^{-u^2/2}), \qquad \quad
\widehat{f_n}(3) \sim \sqrt{2/\pi}(-1 +
2(1-u^2)e^{-u^2/2})n^{-3/2}.
\]
But the weight of $f_n$ at level 1 is $n \cdot \widehat{f_n}(1)^2$
and the weight of $f_n$ at level 3 is ${n \choose 3} \cdot
\widehat{f_n}(3)^2 \sim (n^3/6) \widehat{f_n}(3)^2$; thus the
above imply~(\ref{eqn:formula}) from Theorem~\ref{thm:counter} in
the limit $n \to \infty$ and the proof is complete.
\end{proof}

\end{document}